\documentclass{amsart}
\usepackage{amssymb}

\usepackage{amscd}
\usepackage{thmdefs}

\input{tcilatex}
\theoremstyle{definition}
\theoremstyle{remark}
\numberwithin{equation}{section}

\input tcilatex

\begin{document}
\title[Measure Theory on Graphs]{Measure Theory on Graphs}
\author{Ilwoo Cho}
\address{Saint Ambrose University, Dep. of Math, Davenport, Iowa, U. S. A.}
\email{chowoo@sau.edu}
\thanks{I really appreciate all supports from Saint Ambrose University. I also thank
Prof. T. Anderson and Prof. V. Vega for the useful conversation.}
\date{}
\dedicatory{}
\thanks{}
\maketitle

\begin{abstract}
The main purpose of this paper is to introduce several measures determined
by the given finite directed graph. We will introduce three measures $\mu
_{G},$ $\mu _{G^{\symbol{94}}}$ and $\mu ,$ for the given graph $G.$ Those
three measures are depending on the new algebraic structures, determined by $%
G,$ so-called the diagram set, the reduced diagram set and the free
semigroupoid, respectively. Let $\Bbb{F}^{+}(G)$ be the free semigroupoid of 
$G$ consisting of the empty word $\emptyset ,$ the vertex set $V(G)$ and the
finite path set $FP(G),$ under the admissibility. For each element $w$ in $%
\Bbb{F}^{+}(G),$ we have the graphical image $\delta _{w}$ of $w$ on $\Bbb{R}%
^{2}.$ The set $D(G)$ of such elements is called the diagram set and it is a
quotient semigroupoid of $\Bbb{F}^{+}(G).$ On $D(G),$ we can define a $G$%
-measure $\mu _{G}$ $=$ $d$ $\cup $ $\Delta ,$ where $d$ is a weighted
degree measure on $V(G)$ and $\Delta $ is the diagram length measure on $%
D(G) $ $\setminus $ $V(G).$ Then it is a positive bounded measure. For the
given graph $G,$ we can determine the suadowed graph $G^{\symbol{94}}$ $=$ $G
$ $\cup $ $G^{-1},$ where $G^{-1}$ is the opposite directed graph (or
shadow) of $G.$ Then similar to the previous construction, we can construct
the free semigroupoid $\Bbb{F}^{+}(G^{\symbol{94}})$ of $G^{\symbol{94}}$
and the diagram set $D(G^{\symbol{94}})$ of $G^{\symbol{94}}.$ Define the
reducing relation on $D(G^{\symbol{94}})$ and construct the reduced diagram
set $D_{r}(G^{\symbol{94}}).$ Then it is a groupoid. On it, define the graph
measure $\mu _{G^{\symbol{94}}}$ $=$ $d$ $\cup $ $\Delta ^{r},$ where $d$ is
given as above and $\Delta ^{r}$ is the reduced-diagram-length measure on $%
D_{r}(G^{\symbol{94}})$ $\setminus $ $V(G^{\symbol{94}}).$ Then it is also a
bounded positive measure on $D_{r}(G^{\symbol{94}}).$ Futhermore, the
measure $\mu _{G}$ can be regarded as a restricted measure of $\mu _{G^{%
\symbol{94}}}$ on $\Bbb{F}^{+}(G)$ $\subset $ $\Bbb{F}^{+}(G^{\symbol{94}}).$
In Chapter 3, we observe the graph integral of graph measurable functions
with respect to the graph measure $\mu _{G^{\symbol{94}}}.$ In Chapter 4, we
introduce the restricted measure of the graph measure, so-called subgraph
measures. In Chapter 5, we will extend the graph measure $\mu _{G^{\symbol{94%
}}}$ as a measure $\mu $ defined on the free semigroupoid $\Bbb{F}^{+}(G^{%
\symbol{94}})$ of $G^{\symbol{94}}.$ Different from the graph measure $\mu
_{G^{\symbol{94}}},$ this extended measure $\mu $ is a locally bounded
positive measure. Also, we briefly consider the integration with respect to
this measure.
\end{abstract}

\strut

Let $G$ be a finite directed graph, with its vertex set $V(G)$ and its edge
set $E(G).$ Throughout this paper, we say that a graph is finite if $\left|
V(G)\right| $ $<$ $\infty $ and $\left| E(G)\right| $ $<$ $\infty .$ Let $v$ 
$\in $ $V(G)$ be a vertex. The vertex $v$ has degree $m$ if it has $m$%
-incident edges and denote it by $\deg (v)$ $=$ $m.$ If $v$ is isolated
without the incident edges, then define $\deg (v)$ $=$ $0.$ On the graph $G,$
we can have a finite path $w$ $=$ $e_{1}$ $e_{2}$ ... $e_{k},$ where $e_{1},$
..., $e_{k}$ are admissible directed edges in $E(G)$. Define the length of $%
w $ by $k$ and denote it by $\left| w\right| $ $=$ $k.$ The length of a
finite path is nothing but the cardinality of the admissible edges
constructing the given finite path. Suppose the finite path $w$ has its
initial vertex (or the source) $v_{1}$ and its terminal vertex (or the
range) $v_{2}.$ Then we write $w$ $=$ $v_{1}$ $w$ $v_{2}$ to emphasize the
initial vertex and the terminal vertex of $w.$ Sometimes, we will denote $w$ 
$=$ $v_{1}$ $w$ (or $w$ $=$ $w$ $v_{2}$), for emphasizing the initial vertex
of $w$ (resp. the terminal vertex of $w$). Notice that the finite path $w$
can be expressed as the word in $E(G).$ Define the set $FP(G)$ of all finite
paths. Then

\strut

\begin{center}
$FP(G)=\cup _{n=1}^{\infty }FP_{n}(G),$
\end{center}

where

\begin{center}
$FP_{n}(G)\overset{def}{=}\{w:\left| w\right| =n\}.$
\end{center}

\strut

Let $w_{1}$ and $w_{2}$ be in $V(G)$ $\cup $ $FP(G)$ and assume that the
word $w_{1}$ $w_{2}$ is again contained in $V(G)$ $\cup $ $FP(G).$ Then we
say that $w_{1}$ and $w_{2}$ are admissible (or they are connected under the
ordering $(1,$ $2),$ where $1$ and $2$ are indices of $w_{1}$ and $w_{2}$).
Note that, in general, even though $w_{1}$ and $w_{2}$ are admissible, $%
w_{2} $ and $w_{1}$ are not admissible. However, if $w_{1}$ and $w_{2}$ are
loop finite paths and if $w_{1}$ and $w_{2}$ are admissible, then $w_{2}$
and $w_{1}$ are admissible, too. However, $w_{1}$ $w_{2}$ $\neq $ $w_{2}$ $%
w_{1},$ in this case, because the directions of them are different and hence
they are different finite path. We say that edges $e_{1},$ ..., $e_{k}$ are
admissible if there exists a finite path $w$ $=$ $e_{1}$ ... $e_{k}$ with $%
\left| w\right| $ $=$ $k.$

\strut

Define the free semigroupoid $\Bbb{F}^{+}(G)$ $=$ $\{\emptyset \}$ $\cup $ $%
V(G)$ $\cup $ $FP(G),$ with its binary operation ($\cdot $), the
admissibility, where $\emptyset $ is the empty word in $V(G)$ $\cup $ $%
FP(G). $ The binary operation is the admissibility (or connectedness) of
elements in $V(G)$ $\cup $ $FP(G)$ defined by

\strut

\begin{center}
$(w_{1},w_{2})\longmapsto \left\{ 
\begin{array}{ll}
w_{1}w_{2} & \text{if }w_{1}w_{2}\in V(G)\cup FP(G) \\ 
\emptyset & \text{otherwise.}
\end{array}
\right. $
\end{center}

\strut \strut \strut

The main purpose of this paper is to introduce the measures on a finite
directed graph. Define the diagram map $\delta $ $:$ $\Bbb{F}^{+}(G)$ $%
\rightarrow $ $\Bbb{F}^{+}(G)$ by $w$ $\mapsto $ $\delta _{w},$ where the
diagram $\delta _{w}$ of $w$ is nothing but the graphical image of $w$ in $%
\Bbb{R}^{2}$, on the given graph $G.$ The image $D(G)$ of $\delta $ is a
subset and it also has its admissibility, as binary operation, inherited by
that of $\Bbb{F}^{+}(G)$. i.e., $\delta _{w}$ $\delta _{w^{\prime }}$ $=$ $%
\delta _{ww^{\prime }}.$ Thus the diagram set $D(G)$ with admissibility is
again an algebraic structure.\strut

\strut \strut

In Chapter 1, we will define the measure $\mu _{G}$ on $D(G).$ It is defined
by $d$ $\cup $ $\Delta ,$ where $d$ is the degree measure on $V(G)$ and $%
\Delta $ is the diagram length measure on $D(G)$ $\setminus $ $V(G).$
Simply, $d$ is determined by degrees of vertices and $\Delta $ is determined
by so-called the diagram length. i.e.,

\strut

\begin{center}
$d(S_{V})=\underset{v\in S_{V}}{\sum }\frac{\deg (v)}{\left| V(G)\right| }$
\ and \ $\Delta (S_{FP})=\underset{w\in S_{FP}}{\sum }L(w),$
\end{center}

\strut

where $\frac{\deg (v)}{\left| V(G)\right| }=\frac{\deg _{in}(v)+\deg
_{out}(v)}{\left| V(G)\right| }$ and $L(w)$ $=$ $W(w)$ $\cdot $ $f(w),$
where $w$ $=$ $\delta _{w}$ $\in $ $D(G)$ and $W$ is the corresponding
weighting map of $w$ and $f$ is the length function on $D_{FP}(G)$ $\overset{%
def}{=}$ $D(G)$ $\setminus $ $V(G).$ Here, the sets $S_{V}$ and $S_{FP}$ are
subsets of $V(G)$ and $D_{FP}(G),$ respectively. So, the bounded positive
measure $\mu _{G}$ is defined by

\strut

\begin{center}
$\mu _{G}\left( S\right) =d(S^{\prime })+\Delta (S^{\prime \prime }),$ for
all $S\subseteq D(G)$
\end{center}

\strut

where $S^{\prime }=V(G)\cap S$ and $S^{\prime \prime }=FP(G)\cap S.$ It is
automatically defined that $\mu _{G}(\emptyset )$ $=$ $0.$

\strut

For the given graph $G,$ we can define the shadow $G^{-1}$ of it, as an
opposite directed graph of $G.$ Also, we can have a new finite directed
graph $G^{\symbol{94}}$ with its vertex set $V(G^{\symbol{94}})$ $=$ $V(G)$
and its edge set $E(G^{\symbol{94}})$ $=$ $E(G)$ $\cup $ $E(G^{-1}).$ It is
called the shadowed graph of $G$. By regarding $G^{\symbol{94}}$ as a finite
directed graph, we can construct the free semigroupoid $\Bbb{F}^{+}(G^{%
\symbol{94}})$ of $G^{\symbol{94}}$ and similar to the previous discussion,
we can get the graph measure $\mu _{G^{\symbol{94}}}$ $=$ $d$ $\cup $ $%
\Delta ^{\symbol{94}}$ on $D_{r}(G^{\symbol{94}}),$ where $D_{r}(G^{\symbol{%
94}})$ is the reduced diagram set. Remark that the free semigroupoid $\Bbb{F}%
^{+}(G^{\symbol{94}})$ of the shadowed graph $G^{\symbol{94}}$ is much
bigger than $\Bbb{F}^{+}(G)$ $\cup $ $\Bbb{F}^{+}(G^{-1}).$ However, the
graph measure $\mu _{G}$ (or $\mu _{G^{-1}}$) is the restricted measure on $%
D(G)$ (resp. $D(G^{-1})$) of $\mu _{G^{\symbol{94}}}$ on $D_{r}(G^{\symbol{94%
}}).$ Here, the reduced diagram set $D_{r}(G^{\symbol{94}})$ is the diagram
set $D(G^{\symbol{94}})$ with the reducing property making all loop finite
paths $\delta _{w^{-1}w}$ and $\delta _{ww^{-1}}$ be the corresponding
initial-and-terminal vertices of them. In Chapter 2, we will consider this
shadowed graph measure $\mu _{G^{\symbol{94}}}$ of $G.$ This measure is also
a bounded positive measure. We can show that if the finite directed graphs $%
G_{1}$ and $G_{2}$ are graph-isomorphic, then the measure spaces $(G_{1}^{%
\symbol{94}},$ $\mu _{G_{1}^{\symbol{94}}})$ and $(G_{2}^{\symbol{94}},$ $%
\mu _{G_{2}^{\symbol{94}}})$ are equivalent.

\strut

In Chapter 3, we will construct the integrals on the given graph measure
space $(G^{\symbol{94}},$ $\mu _{G^{\symbol{94}}}).$ In Chapter 4, we
observe the restricted measures of $\mu _{G^{\symbol{94}}},$ induced by the
subgraphs of $G^{\symbol{94}}.$

\strut

Finally, in Chapter 5, we construct a new graph measure $\mu $ $=$ $d$ $\cup 
$ $\omega $ now on the free semigroupoid $\Bbb{F}^{+}(G^{\symbol{94}})$ of
the shadowed graph $G^{\symbol{94}}.$ The construction is very similar to $%
\mu _{G^{\symbol{94}}}.$ However, this measure $\mu $ on $\Bbb{F}^{+}(G^{%
\symbol{94}})$ is just a locally bounded positive measure, in general. Here $%
d$ is the degree measure on $V(G^{\symbol{94}})$ defined before. The
weighted length measure $\omega $ is defined similar to $\Delta ^{r}$ as
follows;

\strut

\begin{center}
$\omega (S)=\underset{w\in S}{\sum }W(\delta _{w}^{r})f(\delta _{w}^{r}),$
for all $S\subseteq FP(G^{\symbol{94}}).$
\end{center}

\strut \strut

\strut

\strut

\section{The Diagram Graph Measures}

\strut

\strut

In this chapter, we will construct the diagram graph measure. Throughout
this chapter, let $G$ be a finite directed graph. The free semigroupoid $%
\Bbb{F}^{+}(G)$ of $G$ is defined as in Introduction. We will determine the
measure $\mu _{G}$ on $D(G),$ where $D(G)$ is the diagram set of $G.$ The
measure space $(D(G),$ $P(D(G))$ $\mu _{G})$ is the main object of this
chapter, where $P(X)$ means the power set of an arbitrary set $X.$

\strut

\strut Let $G$ be a finite directed graph, with its vertex set $V(G)$ and
its edge set $E(G).$ Let $v$ $\in $ $V(G)$ be a vertex. We say that the
vertex $v$ has degree $m$ if it has $m$-incident edges and denote it by $%
\deg (v)$ $=$ $m.$ There are two kinds of incident edges of $v,$ depending
on the direction of $G.$ First, there are inner-incident edges $e$ such that 
$e=e$ $v.$ Second, there are outer-incident edges $e^{\prime }$ such that $%
e^{\prime }$ $=$ $v$ $e^{\prime }.$ The inner-degree $\deg _{in}(v)$ of the
vertex $v$ is defined by

$\strut $

\begin{center}
$\deg _{in}(v)$ $=$ $\left| \{e\in E(G):e=ev\}\right| $
\end{center}

\strut

and the outer-degree $\deg _{out}(v)$ of $v$ is defined by

\strut

\begin{center}
$\deg _{out}(v)=\left| \{e\in E(G):e=ve\}\right| .$
\end{center}

\strut

The degree $\deg (v)$ of $v$ is the sum of $\deg _{in}(v)$ and $\deg
_{out}(v),$

\strut

\begin{center}
$\deg (v)=\deg _{in}(v)+\deg _{out}(v).$
\end{center}

\strut \strut

If $v$ is isolated without the incident edges, then $\deg (v)$ $=$ $0.$
Since the given graph $G$ is finite,

\strut

\begin{center}
$\deg (v)<\infty ,$ \ \ for all \ \ $v\in V(G)$
\end{center}

and

\begin{center}
$\underset{v\in V(G)}{\sum }\frac{\deg (v)}{\left| V(G)\right| }\leq \left|
V(G)\right| \cdot \left| E(G)\right| <\infty .$\ 
\end{center}

\strut \strut \strut

\strut

\begin{definition}
Let $G$ be the given finite directed graph and $\Bbb{F}^{+}(G),$ the
corresponding free semigroupoid of $G$ and let $w$ $\in $ $\Bbb{F}^{+}(G)$.
Define the diagram $\delta _{w}$ of $w$ by the graphical image of $w$ on the
graph $G.$ We say that the element $w$ in $\Bbb{F}^{+}(G)$ is basic if the
diagram $\delta _{w}$ of it is $w,$ itself, on $G.$ Define the set $D_{k}(G)$
of all diagrams with length $k$ by

\strut

$\ \ \ \ \ \ \ \ \ \ \ \ \ \ \ \ \ \ \ D_{k}(G)=\{\delta _{w}:\left| \delta
_{w}\right| =k\},$ \ \ for all \ $k\in \Bbb{N}.$

\strut

Remark that, by regarding $\delta _{w}$ as a finite path, we can determine
the length $\left| \delta _{w}\right| $ of the diagram $\delta _{w}.$ The
set $D(G)$ of all diagrams is defined by $D(G)$ $=$ $\{\emptyset \}$ $\cup $ 
$(\cup _{n=0}^{\infty }$ $D_{n}(G)),$ where

\strut

$\ \ \ \ \ \ \ \ \ \ \ \ \ \ \ \ \ \ \ \ \ D_{0}(G)$ $\overset{def}{=}$ $%
\{\delta _{v}$ $:$ $v$ $\in $ $V(G)\}=V(G).$

\strut

Denote $D(G)$ $\setminus $ $\left( V(G)\cup \{\emptyset \}\right) $ by $%
D_{FP}(G).$ i.e., $D_{FP}(G)$ $=$ $D(G)$ $\cap $ $FP(G).$
\end{definition}

\strut

Notice that there exists $N$ $\in $ $\Bbb{N}$ such that

\strut

\begin{center}
$D(G)=\cup _{n=0}^{N}D_{n}(G),$
\end{center}

\strut

since $G$ is a finite graph.

\strut

Let's discuss about the differences between $FP(G)$ and $D_{FP}(G).$ Assume
that $l$ $\in $ $E(G)$ is a loop-edge. Then it has its diagram $\delta _{l}$ 
$=$ $l$ on $G,$ because the diagram $\delta _{l}$ of $l$ is $l,$ itself on
the given graph $G.$ In fact, by definition, $\delta _{e}$ $=$ $e,$ for all $%
e$ $\in $ $E(G).$ In other words, $E(G)$ $=$ $D_{1}(G)$ $\subset $ $D(G).$
So, as a finite path, the loop-edge $l$ is basic. However, the diagram $%
\delta _{l^{2}}$ of the finite path $l^{2}$ (with its length $2$) is $l$ $=$ 
$\delta _{l}$ on $G,$ because the finite path $l^{2}$ also has its graphical
image $l$ on $G.$ So, in this case, the diagram $\delta _{l^{2}}$ of $l^{2}$
is exactly same as the diagram $\delta _{l}$ of the basic finite path $l.$
More generally, all diagrams $\delta _{l^{k}}$ of loop\ finite paths $l^{k}$%
, $k$ $\in $ $\Bbb{N}$ $\setminus $ $\{1\},$ are same as $l$ $=$ $\delta
_{l} $ of the basic finite path $l.$

\strut

\begin{example}
Let $e=v_{1}ev_{2}\in E(G),$ with $v_{1}\neq v_{2}$ in $V(G)$ and let $l$ $=$
$v_{2}lv_{2}$ $\in $ $FP(G)$ be a loop finite path. Then the finite paths $e$
and $l$ are admissible. i.e., the finite path $el$ $\in $ $FP(G).$ Assume
that there exists a basic loop finite path $w$ $=$ $v_{2}$ $w$ $v_{2}$ such
that $l$ $=$ $w^{k},$ for some $k$ $\in $ $\Bbb{N},$ in $FP(G).$ Then the
diagram $\delta _{el}$ of $el$ is same as the diagram $\delta _{ew},$ and
the diagram $\delta _{ew}$ is identified with $ew$ on $G$.\strut
\end{example}

\strut

Assume that $v\in \Bbb{F}^{+}(G)$ is a vertex of $G.$ Then the diagram $%
\delta _{v}$ of $v$ is always $v$ on $G.$ Also, if $e$ $\in $ $\Bbb{F}%
^{+}(G) $ is an edge contained in $E(G),$ then the diagram $\delta _{e}$ of $%
e$ is $e $ on $G.$ Now, we can define the diagram map

$\strut $

\begin{center}
$\delta $ $:$ $\Bbb{F}^{+}(G)$ $\rightarrow $ $D(G)$ by the map sending $w$
to $\delta _{w},$
\end{center}

\strut

which is a surjection from $\Bbb{F}^{+}(G)$ onto $D(G).$ If $w_{1}$ and $%
w_{2}$ have the same diagram on $G,$ then $\delta _{w_{1}}$ $=$ $\delta
_{w_{2}}$ on $G.$ We can easily get that;

\strut

\begin{lemma}
(1) If $v\in V(G),$ then $\delta _{v}=v$ on $G.$

(2) If $e\in E(G),$ then $\delta _{e}=e$ on $G.$

(3) An element $w$ $\in $ $\Bbb{F}^{+}(G)$ is basic if and only if $\delta
_{w}=w$ on $G.$

(4) If $w\in \Bbb{F}^{+}(G)$ is not basic, then there exists a unique basic
finite path $w_{0}$ $\in $ $\Bbb{F}^{+}(G)$ such that $\delta _{w}$ $=$ $%
\delta _{w_{0}}$ on $G.$
\end{lemma}

\strut

\begin{proof}
It suffices to prove the uniqueness in (4). The existence is guaranteed by
(3). i.e., the diagram $\delta _{w}$ of $w$ in $D(G)$ is the basic element
having the same diagram with that of $w.$ i.e., $\delta _{w}$ $=$ $\delta
_{\delta _{w}}.$ Suppose that $w_{1}$ and $w_{2}$ are basic elements such
that $\delta _{w_{1}}$ $=$ $\delta _{w}$ $=$ $\delta _{w_{2}}.$ By the
basicness of $w_{1}$ and $w_{2},$ we have $w_{1}$ $=$ $\delta _{w_{1}}$ $=$ $%
\delta _{w_{2}}$ $=$ $w_{2}.$
\end{proof}

\strut \strut \strut

By definition of finite directed trees, $D(G)$ $=$ $\Bbb{F}^{+}(G),$
whenever $G$ is a tree.

\strut

\begin{proposition}
Let $G$ be a finite directed tree. Then $\Bbb{F}^{+}(G)=D(G).$ $\square $
\end{proposition}

\strut \strut \strut

We will construct the measure $\mu _{G}$ on the free semigroupoid $D(G)$ of
the finite directed graph $G.$ First, consider the length function $f$ on $%
D(G).$ Define this map

\strut

\begin{center}
$f:D(G)\rightarrow \Bbb{R}$
\end{center}

by

\begin{center}
$f(w)=k,$ \ for all \ $w\in D_{k}(G).$
\end{center}

\strut

Remark that there exists $N$ $\in $ $\Bbb{N}$ such that $D(G)$ $=$ $\cup
_{n=0}^{N}$ $D_{n}(G)$, for any finite graph $G.$ Moreover, since $G$ is
finite $\left| D_{k}(G)\right| $ $<$ $\infty ,$ for all $k$ $=$ $0,$ $1,$
..., $N.$ So, we can define the following map $F$ $:$ $P\left( D(G)\right) $ 
$\rightarrow $ $\Bbb{N},$

\strut

\begin{center}
$F\left( S\right) =\underset{\delta \in S}{\sum }f(\delta ),$
\end{center}

\strut

for all $S\in P\left( D(G)\right) ,$ where $P\left( D(G)\right) $ is the
power set of $D(G).$ For instance, if $S$ $=$ $\{\delta _{1},$ $\delta
_{2}\} $ and if $\delta _{j}$ $\in $ $D_{k_{j}}(G),$ for $j$ $=$ $1,$ $2,$
then

\strut

\begin{center}
$F(S)=k_{1}+k_{2}\in \Bbb{N}.$
\end{center}

\strut

This map $F$ is called the diagram-length map.

\strut \strut

\begin{lemma}
Let $G$ be a finite directed graph and let $D(G)$ be the set of all diagrams
on $G.$ Then the map $F$ is bounded. i.e., $F(S)$ $<$ $\infty ,$ for all $S$ 
$\in $ $P\left( D(G)\right) .$
\end{lemma}

\strut

\begin{proof}
Let's take $S$ $=$ $D(G),$ the largest subset of $D(G).$ It is sufficient to
show that $F$ $\left( D(G)\right) $ $<$ $\infty .$ Observe

\strut

$\ \ \ \ \ \ F\left( D(G)\right) =F\left( \cup _{k=0}^{N}D_{k}(G)\right)
=\sum_{k=0}^{N}F\left( D_{k}(G)\right) $

\strut

for some $N\in \Bbb{N},$ since $D_{k}(G)$'s are mutually disjoint

\strut

$\ \ \ \ \ \ \ \ \ \ \ \ \ \ \ \ \ \ \ \ \ \ =0\cdot \left| V(G)\right|
+\sum_{k=1}^{N}k\left| D_{k}(G)\right| $

\strut

since $F\left( D_{k}(G)\right) =\underset{\delta \in D_{k}(G)}{\sum }%
f(\delta )=k\cdot \left| D_{k}(G)\right| ,$ for all $k$ $=$ $0,$ $1,$ ..., $%
N,$

\strut

$\ \ \ \ \ \ \ \ \ \ \ \ \ \ \ \ \ \ \ \ \ \ <\infty ,$

\strut

since $G$ is finite and hence $D_{k}(G)$ has the bounded cardinality, for
all $k$ $=$ $1,$ ..., $N.$
\end{proof}

\strut

Define the weighted diagram-length function $L$ on $D(G).$

\strut \strut

\begin{definition}
Let $G$ be a finite directed graph and $\Bbb{F}^{+}(G)$, the corresponding
free semigroupoid of $G$ and let $D(G)$ be the set of all diagrams on $G.$
Let $F$ $:$ $P(D_{FP}(G))$ $\rightarrow $ $\Bbb{N}$ be the diagram-length
map defined above. Let $w$ $=$ $e_{1}$ ... $e_{k}$ be a $\emph{basic}$
finite path with admissible edges $e_{1},$ ..., $e_{k},$ where $k$ $\in $ $%
\Bbb{N}.$ Define the weighting $W_{E}$ $:$ $E(G)$ $\rightarrow $ $(0,$ $1)$
by $W_{E}(e)$ $=$ $r_{e}$, where $(0,$ $1)$ is the open interval in $\Bbb{R}%
. $ More generally, define the weighting map $W$ $:$ $D_{FP}(G)$ $%
\rightarrow $ $(0,$ $1)$ by

\strut

\ $\ \ \ \ \ \ \ \ \ W(w)$ $=$ $\left\{ 
\begin{array}{ll}
W_{E}(w) & \text{if }w\in E(G) \\ 
&  \\ 
\Pi _{j=1}^{k}W_{E}(e_{j}) & \text{if }w=e_{1}...e_{k}\in D_{FP}(G)
\end{array}
\right. $

\strut \strut

Define the weighted diagram-length map (in short, WDL-map) $L$ $:$ $P\left(
D_{FP}(G)\right) $ $\rightarrow $ $(0,$ $1)$ by

\strut

$\ \ \ \ \ \ \ \ \ \ \ \ \ \ L(S)\overset{def}{=}\underset{w\in S}{\sum }%
W(w)f(w),$ for all $S\in P\left( D_{FP}(G)\right) .$

\strut

The unweighted diagram-length map (in short, UWDL-map) $L_{u}$ on $%
P(D_{FP}(G))$ is defined by $L_{u}$ $(S)$ $=$ $F$ $\left( S\right) ,$ for
all $S$ in $P\left( D_{FP}(G)\right) .$ \strut
\end{definition}

\strut \strut \strut \strut \strut

By definition, we can get that;

\strut

\begin{lemma}
Let $G$ be a finite directed graph and the map $L$ $:$ $P\left(
D_{FP}(G)\right) $ $\rightarrow $ $\Bbb{R}$ be the WDL-map. Then it is
bounded. The UWDL-map $L_{u}$ is also bounded.
\end{lemma}

\strut

\begin{proof}
It suffices to show that $L\left( D_{FP}(G)\right) $ is bounded. By
definition, we have that

\strut

$\ \ \ \ \ \ \ \ \ \ \ \ \ \ \ \ \ \ \ L\left( D_{FP}(G)\right) =\underset{%
w\in D_{FP}(G)}{\sum }W(w)f(w)<\infty ,$

\strut

since the diagram set $D(G)$ has a finite cardinality and so does $D_{FP}(G)$%
. Similarly, we can show that

\strut

$\ \ \ L_{u}$ $\left( D_{FP}(G)\right) $ $=$ $F$ $\left( D_{FP}(G))\right) $ 
$=$ $F\left( \cup _{k=1}^{N}D_{k}(G)\right) <\infty .$

\strut
\end{proof}

\strut \strut \strut \strut \strut \strut

Now, we will define the measure $\mu _{G}$ on $D(G).$

\strut \strut

\begin{definition}
Let $G$ be a finite directed graph and let $D(G)$ be the corresponding
diagram set of $G.$ Define the degree measure $d$ on the vertex set $V(G)$ by

\strut

$\ \ \ \ \ \ \ \ \ \ \ \ \ \ \ \ \ \ d(S)=\underset{v\in S}{\sum }\frac{\deg
(v)}{\left| V(G)\right| },$ \ \ for all \ \ $S\subseteq V(G).$

\strut

Also, define the weighted measure $\Delta $ on the finite path set $D(G)$ $%
\cap $ $FP(G)$ by

\strut

\ \ \ \ \ \ \ \ \ $\ \ \ \ \ \ \ \ \ \ \Delta (S)=L(S),$ \ for all \ $%
S\subseteq D_{FP}(G).$

\strut

The unweighted measure $\Delta _{u}$ on $FP(G)$ by

\strut

\ \ \ \ \ \ \ \ \ \ \ $\ \ \ \ \ \ \ \Delta _{u}(S)=L_{u}(S),$ \ for all \ $%
S\subseteq D_{FP}(G).$

\strut

Finally, we have the (weighted) graph measure (or $G$-measure) $\mu _{G}$ on 
$D(G),$ denoted by $d$ $\cup $ $\Delta $, defined by

\strut

$\ \ \ \ \ \ \ \ \ \ \ \ \ \mu _{G}(S)=d(S_{V})+\Delta (S_{FP}),$ for all $%
S\subseteq D(G),$

\strut

where $S_{V}=S\cap V(G)$ and $S_{FP}=S\cap D_{FP}(G).$ Similarly, the
unweighted $G$-measure $\mu _{G}^{u}$ is defined by $d$ $\cup $ $\Delta
_{u}. $
\end{definition}

\strut

The $G$-measure $\mu _{G}$ and $\mu _{G}^{u}$ are indeed bounded measures on 
$D(G).$ Instead of using $P(D(G)),$ we will use the quotient set $P\left( 
\Bbb{F}^{+}(G)\right) $ $/$ $\mathcal{R}.$

\strut

\begin{theorem}
The $G$-measure $\mu _{G}$ is a bounded positive measure on $P\left(
D(G)\right) .$ The unweighted $G$-measure $\mu _{G}^{u}$ is also a bounded
positive measure on $P\left( D(G)\right) .$
\end{theorem}

\strut

\begin{proof}
To show that $\mu _{G}=d\cup \Delta $ is a measure on $D(G),$ we will check
the followings;

\strut

(i) \ \ $\mu _{G}(S)\geq 0,$ for all $S\subseteq D(G).$ If $S$ is a subset
of $D(G),$ then there always is a separation $S_{V}$ and $S_{FP}$ of $S$,
where $S_{V}$ $=$ $V(G)$ $\cap $ $S$ and $S_{FP}$ $=$ $D_{FP}(G)$ $\cap $ $S$%
, such that $S_{V}$ $\cup $ $S_{FP}$ $=$ $S$ and $S_{V}$ $\cap $ $S_{FP}$ $=$
$\varnothing .$ So,

\strut

$\ \ \ \ \ \ \ \mu _{G}(S)$ $=d(S_{V})+\Delta (S_{FP})=\underset{v\in S_{V}}{%
\sum }\frac{\deg (v)}{\left| V(G)\right| }+\underset{w\in S_{FP}}{\sum }%
W(w)f(w)\geq 0.$

\strut

(ii) \ For any $S\subseteq D(G),$ $\mu _{G}(S)<\infty $. We need to show
that $\mu _{G}\left( D(G)\right) $ $<$ $\infty .$ We have that

\strut

$\ \ \ \ \ \ \ \ \ \ \ \ \ \ \ \ \ \ \mu _{G}\left( D(G)\right) =d\left(
V(G)\right) +\Delta \left( D_{FP}(G)\right) .$

\strut

By the previous lemma, $\Delta \left( D_{FP}(G)\right) <\infty $. Also,

\strut

$\ \ \ \ \ \ \ \ \ \ \ d\left( V(G)\right) =$ $\underset{v\in V(G)}{\sum }%
\frac{\deg (v)}{\left| V(G)\right| }\leq \left| V(G)\right| \cdot \left|
E(G)\right| <\infty ,$

\strut

since $G$ is a finite graph. Therefore, $\mu _{G}\left( D(G)\right) <\infty
. $

\strut

(iii) Let $(S_{n})_{n=1}^{\infty }$ be a mutually disjoint sequence of
subsets of $D(G).$ Then

\strut

$\ \ \ \mu _{G}\left( \cup _{n=1}^{\infty }S_{n}\right) =\mu _{G}\left(
\left( \cup _{n=1}^{\infty }S_{n,V}\right) \cup \left( \cup _{n=1}^{\infty
}S_{n,FP}\right) \right) $

\strut

where $S_{n,V}=V(G)\cap S_{n}$ and $S_{n,FP}=D_{FP}(G)\cap S_{n},$ for all $%
n $ $\in $ $\Bbb{N}$

\strut

$\ \ \ \ \ \ \ \ \ \ \ \ \ =d\left( \cup _{n=1}^{\infty }S_{n,V}\right)
+\Delta \left( \cup _{n=1}^{\infty }S_{n,FP}\right) $

\strut

by the definition of $\mu _{G}=d\cup \Delta $

\strut

$\ \ \ \ \ \ \ \ \ \ \ \ \ =\sum_{n=1}^{\infty }d\left( S_{n,V}\right)
+\sum_{n=1}^{\infty }\Delta \left( S_{n,FP}\right) $

\strut

by the definition of $d$ and $\Delta $ and by the disjointness of $%
\{S_{n,V}\}_{n=1}^{\infty }$ and $\{S_{n,FP}\}_{n=1}^{\infty }$

\strut

$\ \ \ \ \ \ \ \ \ \ \ \ \ =\sum_{n=1}^{\infty }\mu _{G}(S_{n}).$

\strut \strut

By (i) and (iii), the $G$-measure $\mu _{G}=d\cup \Delta $ is a measure on $%
D(G),$ and by (ii), it is a bounded measure on $D(G).$
\end{proof}

\strut \strut \strut

\begin{definition}
Let $G$ be a finite directed graph. The triple $(D(G),$ $P\left( D(G)\right)
,$ $\mu _{G})$ is called the graph measure space (in short, $G$-measure
space). For convenience, the triple is denoted by $(G,$ $\mu _{G}).$
\end{definition}

\strut \strut \strut

\begin{proposition}
Let $(G,$ $\mu _{G})$ be a $G$-measure space. Then the $G$-measure $\mu _{G}$
is atomic, in the sense that, for all $S$ $\subseteq $ $D(G)$ $\setminus $ $%
\{\emptyset \},$ $0$ $<$ $\mu _{G}(S)$ $<$ $\infty .$ $\square $
\end{proposition}

\strut \strut \strut \strut

\strut \strut

\strut \strut

\section{Reduced Diagram Graph Measures}

\strut

\strut

Throughout this chapter, let $G$ be a finite directed graph and let $\Bbb{F}%
^{+}(G)$ be the corresponding free semigroupoid. In this chapter, we will
construct the shadowed graph $G^{\symbol{94}}$ of the given graph $G,$ as a
new finite directed graph. By constructing the free semigroupoid $\Bbb{F}%
^{+}(G^{\symbol{94}})$ of the shadowed graph $G^{\symbol{94}},$ we can
define the free groupoid $\Bbb{F}^{+}(G^{\symbol{94}})$ of the given graph $%
G.$ Similar to the previous chapter, we can have the $G^{\symbol{94}}$%
-measure space

\strut

\begin{center}
$\left( G^{\symbol{94}},\text{ }\mu _{G^{\symbol{94}}}\right) =\left(
D_{r}(G^{\symbol{94}}),\text{ }P\left( D_{r}(G^{\symbol{94}})\right) ,\text{ 
}\mu _{G^{\symbol{94}}}\right) .$
\end{center}

\strut

But the graph measure $\mu _{G^{\symbol{94}}}$ is defined slightly
differently. To define this $G^{\symbol{94}}$-measure $\mu _{G^{\symbol{94}%
}} $ $=$ $d$ $\cup $ $\Delta ^{\symbol{94}}$, (in particular, to define $%
\Delta ^{\symbol{94}}$), we will use the reduced diagrams $D_{r}(G^{\symbol{%
94}}),$ instead of using the diagrams $D(G^{\symbol{94}}).$

\strut

\strut

\strut

\subsection{Shadowed Graphs}

\strut

\strut

Let $G$ be the given finite directed graph and $\Bbb{F}^{+}(G),$ the
corresponding free semigroupoid.

\strut

\begin{definition}
Let $G$ be the given graph. Then the shadow $G^{-1}$ is defined by a graph
with

\strut

$\ \ \ \ \ \ \ \ \ V(G^{-1})=V(G)$ \ and \ $E(G^{-1})=\{e^{-1}:e\in E(G)\},$

\strut

where $e^{-1}$ is the opposite directed edge of $e.$ i.e., the shadow $%
G^{-1} $ is the opposite directed graph of $G.$
\end{definition}

\strut

Let $e=v_{1}ev_{2},$ with $v_{1},$ $v_{2}$ $\in $ $V(G).$ Then the opposite
directed edge $e^{-1}$ of $e$ satisfies that $e^{-1}$ $=$ $v_{2}$ $e^{-1}$ $%
v_{1}.$ Notice that the admissibility on $G^{-1}$ is preserved oppositely by
that of $G.$ i.e., the finite path $w^{-1}$ $=$ $e_{1}^{-1}$ ... $e_{k}^{-1}$
is in $FP(G^{-1})$ if and only if $w$ $=$ $e_{k}$ ... $e_{1}$ is in $FP(G).$
More generally, the element $w^{-1}$ $=$ $w_{1}^{-1}$ ... $w_{k}^{-1}$ is in 
$\Bbb{F}^{+}(G^{-1})$ if and only if $w$ $=$ $w_{k}$ ... $w_{1}$ is in $\Bbb{%
F}^{+}(G).$

\strut

\begin{proposition}
$\left| \Bbb{F}^{+}(G^{-1})\right| =\left| \Bbb{F}^{+}(G)\right| .$
\end{proposition}

\strut

\begin{proof}
By definition, the free semigroupid $\Bbb{F}^{+}(G^{-1})$ of $G^{-1}$ is $%
\{\emptyset \}$ $\cup $ $V(G)$ $\cup $ $FP(G^{-1}).$ So, it suffices to show
that

\strut

$\ \ \ \ \ \ \ \ \ \ \ \ \ \ \ \ \ \ \ \ \ \ \ \left| FP(G^{-1})\right|
=\left| FP(G)\right| .$

\strut

There exists a natural map from $FP(G)$ onto $FP(G^{-1})$ defined by $w$ $%
\mapsto $ $w^{-1}.$ Since the admissibility on $G^{-1}$ is preserved
oppositely by that of $G$ and since we can regard $G$ as $(G^{-1})^{-1},$
the map $w$ $\mapsto $ $w^{-1}$ is bijective. So, the above equality holds.
\end{proof}

\strut

Since $\left| \delta \left( \Bbb{F}^{+}(G^{-1})\right) \right| =\left|
\delta \left( \Bbb{F}^{+}(G)\right) \right| ,$ where $\delta $ is the
diagram map, we can get that;

\strut

\begin{corollary}
$\left| D(G^{-1})\right| =\left| D(G)\right| .$ \ $\square $
\end{corollary}

\strut

By the previous two results, we can get that;

\strut

\begin{theorem}
The graph measure spaces $(G,$ $\mu _{G})$ and $(G^{-1},$ $\mu _{G^{-1}})$
are equivalent. $\square $
\end{theorem}

\strut

Now, let's define the shadowed graph $G^{\symbol{94}}$ of $G.$

\strut

\begin{definition}
Let $G$ be given as above and let $G^{-1}$ be the shadow of $G.$ Define the
shadowed graph $G^{\symbol{94}}$ by the directed graph with

\strut

$\ \ \ \ V(G^{\symbol{94}})=V(G)=V(G^{-1})$ \ and \ $E(G^{\symbol{94}%
})=E(G)\cup E(G^{-1}).$
\end{definition}

\strut

Notice that $\Bbb{F}^{+}(G^{\symbol{94}})\neq \Bbb{F}^{+}(G)\cup \Bbb{F}%
^{+}(G^{-1}),$ since we can also take a mixed finite path $w_{1}^{i_{1}}$ $%
w_{2}^{i_{2}}$ ... $w_{m}^{i_{m}},$ where $(i_{1},$ ..., $i_{m})$ $\in $ $%
\{1,$ $-1\}^{m}$ is a mixed $m$-tuple of $1$ and $-1,$ in $\Bbb{F}^{+}(G^{%
\symbol{94}}).$ Such mixed finite paths are not contained in $\Bbb{F}^{+}(G)$
$\cup $ $\Bbb{F}^{+}(G^{-1}).$ Therefore, in general,

\strut

\begin{center}
$\Bbb{F}^{+}(G^{\symbol{94}})\supsetneqq \left( \Bbb{F}^{+}(G)\cup \Bbb{F}%
^{+}(G^{-1})\right) .$
\end{center}

\strut \strut

\strut

\strut

\subsection{Reduced Diagrams on the Shadowed Graph}

\strut

\strut

Let $G$ be the given finite directed graph and $G^{-1},$ the shadow of $G$,
and let $G^{\symbol{94}}$ be the shadowed graph of $G.$ Like Chapter 1, we
can define the diagram $D(G^{\symbol{94}})$ of $G^{\symbol{94}},$ as the
image of the diagram map $\delta $ $:$ $\Bbb{F}^{+}(G^{\symbol{94}})$ $%
\rightarrow $ $D(G^{\symbol{94}})$ defined by $w$ $\mapsto $ $\delta _{w}.$
On this diagram set $D(G^{\symbol{94}})$, we will give the reducing process
and we will define the reduced diagram set $D_{r}(G^{\symbol{94}})$ of $G^{%
\symbol{94}}.$ This reduced diagram set $D_{r}(G^{\symbol{94}})$ is a
groupoid, under the admissibility inherited by that of $D(G),$ with the
reducing property.

\strut

\begin{definition}
Let $G$, $G^{-1}$ and $G^{\symbol{94}}$ be given as above and let $\Bbb{F}%
^{+}(G^{\symbol{94}})$ be the free semigroupoid of $G^{\symbol{94}}.$ Let $%
\delta $ $:$ $\Bbb{F}^{+}(G^{\symbol{94}})$ $\rightarrow $ $D(G^{\symbol{94}%
})$ be the diagram map defined in Chapter 1, for the finite directed graph $%
G^{\symbol{94}}.$ Define the reduced diagram map $\delta ^{r}$ $:$ $\Bbb{F}%
^{+}(G^{\symbol{94}})$ $\rightarrow $ $D(G^{\symbol{94}})$ by

\strut

\ \ \ $\ \ \ \ \ \ \ \ \ \ \ \ \ \ \ w\longmapsto \delta _{w}^{r}=\delta
_{w},$ $\forall $ $w$ $\in \Bbb{F}^{+}(G)\cup \Bbb{F}^{+}(G^{-1})$

and

\strut

(2.1)$\ \ \ \ \ \ \ \ w^{-1}w\mapsto \delta _{w^{-1}w}^{r}=\delta _{v}$ \
and\ \ $ww^{-1}\mapsto \delta _{w^{-1}w}^{r}=\delta _{v^{\prime }},$

\strut

for all $w=v^{\prime }wv\in FP(G)$ and $w^{-1}=vw^{-1}v^{\prime }\in $ $%
FP(G^{-1}),$ with $v,$ $v^{\prime }$ $\in $ $V(G),$ and

\strut

$\ \ \ \ \ \ \ \ \ \ \ \ \ \ \ \ \ \ \
w_{1}^{i_{1}}w_{2}^{i_{2}}...w_{n}^{i_{n}}\longmapsto \delta
_{w_{1}^{i_{1}}w_{2}^{i_{2}}...w_{n}^{i_{n}}}^{r}=\delta _{\delta
_{w_{1}^{i_{1}}...w_{n}^{i_{n}}}}^{r},$

\strut

for the $n$-tuple $(i_{1},$ ..., $i_{n})$ $\in $ $\{1,$ $-1\}^{n},$ for all $%
n$ $\in $ $\Bbb{N}$. The image $\delta _{w}^{r}$ of $w$ in $\Bbb{F}^{+}(G)$
is called the reduced diagram of $w.$ Define the set $D_{r}(G)$ of all
reduced diagrams by the image $\delta ^{r}\left( \Bbb{F}^{+}(G^{\symbol{94}%
})\right) $ of $\delta ^{r}.$ i.e.,

\strut

$\ \ \ \ \ \ \ \ \ \ \ D_{r}(G)$ $\overset{def}{=}$ $\delta ^{r}\left( \Bbb{F%
}^{+}(G^{\symbol{94}})\right) =\delta ^{r}\left( \delta (\Bbb{F}^{+}(G^{%
\symbol{94}}))\right) =\delta ^{r}\left( D(G^{\symbol{94}})\right) .$
\end{definition}

\strut \strut

\strut

\strut \strut

\subsection{Reduced Diagram Graph Measure Spaces}

\strut

\strut

Let $G,$ $G^{-1}$ and $G^{\symbol{94}}$ be given as above. Then we can
define the reduced diagram graph measure $\mu _{G^{\symbol{94}}}$ by $\mu
_{G^{\symbol{94}}}$ $=$ $d$ $\cup $ $\Delta ^{\symbol{94}}.$

\strut

\begin{definition}
Let $G^{\symbol{94}}$ be the shadowed graph of a finite directed graph $G$
and let $\Delta ^{\symbol{94}}$ be the reduced diagram length measure
defined by the diagram length measure on $D_{r}^{FP}(G^{\symbol{94}})$ $%
\overset{def}{=}$ $D_{r}(G^{\symbol{94}})$ $\cap $ $FP(G^{\symbol{94}})$.
Define the reduced diagram graph measure $\mu _{G^{\symbol{94}}}$ by $d$ $%
\cup $ $\Delta ^{\symbol{94}}.$ i.e.,

\strut

$\ \ \ \ \ \ \ \ \ \ \ \ \ \ \ \ \ \ \ \ \ \mu _{G^{\symbol{94}}}\left(
S\right) =d\left( S_{V}\right) +\Delta ^{\symbol{94}}\left( S_{FP}\right) ,$

\strut

for all $S\subseteq D_{r}(G^{\symbol{94}}),$ where $S_{V}$ $=$ $S$ $\cap $ $%
V(G^{\symbol{94}})$ and $S_{FP}$ $=$ $S$ $\cap $ $D_{r}^{FP}(G^{\symbol{94}%
}).$ The triple $(D_{r}(G^{\symbol{94}}),$ $P$ $\left( D_{r}(G^{\symbol{94}%
})\right) ,$ $\mu _{G^{\symbol{94}}})$ is called the reduced diagram graph
measure space of the graph $G$ (in short, $G^{\symbol{94}}$-measure space of 
$G$), denoted by $(G^{\symbol{94}},$ $\mu _{G^{\symbol{94}}}).$
\end{definition}

\strut

Then, similar to Chapter 1, we can get that;

\strut

\begin{theorem}
Let $(G^{\symbol{94}},$ $\mu _{G^{\symbol{94}}})$ be the $G^{\symbol{94}}$%
-measure space of $G.$ Then $\mu _{G^{\symbol{94}}}$ is bounded positive
measure on $D_{r}(G^{\symbol{94}}).$ In particular, this measure is atomic,
in the sense that, all subsets $S$ in $P$ $\left( D_{r}(G^{\symbol{94}%
})\right) ,$ which is not $\{\emptyset \},$ have bounded measures. $\square $
\end{theorem}

\strut \strut \strut

The key idea to prove the above theorem is to understand the reduced diagram
set $D_{r}(G^{\symbol{94}})$ has finite cardinality. We can use the same
techniques used in Chapter 1.\strut By definition, the following proposition
is easily proved.

\strut

\begin{proposition}
Let $\mu _{G}$ and $\mu _{G^{-1}}$ be the $G$-measure on $D(G)$ and $\,$the $%
G^{-1}$-measure on $D(G^{-1}),$ respectively. Then the $G^{\symbol{94}}$%
-measure $\mu _{G^{\symbol{94}}}$ satisfies that $\mu _{G^{\symbol{94}}}$ $=$
$\mu _{G}$ on $D(G)$ and $\mu _{G^{\symbol{94}}}$ $=$ $\mu _{G^{-1}}$ on $%
D(G^{-1}).$ $\square $
\end{proposition}

\strut \strut

Let $G_{1}$ and $G_{2}$ be finite directed graph and $G_{1}^{\symbol{94}}$
and $G_{2}^{\symbol{94}},$ the shadowed graphs of $G_{1}$ and $G_{2},$
respectively. Assume that $G_{1}$ and $G_{2}$ are graph-isomorphic. i.e.,
there exists a graph-isomorphism $\varphi $ $:$ $G_{1}$ $\rightarrow $ $%
G_{2} $ such that (i) $\varphi $ is the bijection from $V(G_{1})$ onto $%
V(G_{2}),$ (ii) $\varphi $ is the bijection from $E(G_{1})$ onto $E(G_{2})$,
(iii) for $v,$ $v^{\prime }$ $\in $ $V(G_{1})$ such that $e$ $=$ $v$ $e$ $%
v^{\prime }$ in $E(G_{1}),$ $\varphi (e)$ $=$ $\varphi (v)$ $\varphi (e)$ $%
\varphi (v^{\prime })$. i.e., $\varphi $ preserves the admissibility. Then
the $G_{1}^{\symbol{94}}$-measure (resp. $G_{1}$-measure) and $G_{2}^{%
\symbol{94}} $-measure (resp. $G_{2}$-measure) are equivalent.

\strut

\begin{theorem}
Let $G_{1}$ and $G_{2}$ be graph-isomorphic finite directed graphs. Then the
shadowed graph measure spaces $(G_{1}^{\symbol{94}},$ $\mu _{G_{1}^{\symbol{%
94}}})$ and $(G_{2}^{\symbol{94}},$ $\mu _{G_{2}^{\symbol{94}}})$ are an
equivalent measure spaces. In particular, the graph measure spaces $(G_{1},$ 
$\mu _{G_{1}})$ and $(G_{2},$ $\mu _{G_{2}})$ are also equivalent measure
spaces. $\square $
\end{theorem}

\strut

By the definition of the graph-isomorphisms, we can easily prove the above
theorem.\strut The main idea is that $D_{r}(G_{1}^{\symbol{94}})$ (resp. $%
D(G_{1})$) and $D_{r}(G_{2}^{\symbol{94}})$ (resp. $D(G_{2})$) are
reduced-diagram-isomorphic (resp. diagram-isomorphic). In fact, we can have
that

\strut

$\ G_{1}$ and $G_{2}$ are graph-isomorphic

\strut

$\ \ \ \Longrightarrow $ $G_{1}^{\symbol{94}}$ and $G_{2}^{\symbol{94}}$ are
graph-isomorphic

\strut

$\ \ \ \Longrightarrow \Bbb{F}^{+}(G_{1}^{\symbol{94}})$ and $\Bbb{F}%
^{+}(G_{2}^{\symbol{94}})$ are free--semigroupoid-isomorphic

\strut

$\ \ \ \Longrightarrow D(G_{1}^{\symbol{94}})$ and $D(G_{2}^{\symbol{94}})$
are diagram-isomorphic

\strut

$\ \ \ \Longrightarrow D_{r}(G_{1}^{\symbol{94}})$ and $D_{r}(G_{2}^{\symbol{%
94}})$ are reduced-diagram-isomorphic

\strut

$\ \ \ \Longrightarrow (G_{1}^{\symbol{94}},$ $\mu _{G_{1}^{\symbol{94}}})$
and $(G_{2}^{\symbol{94}},$ $\mu _{G_{2}^{\symbol{94}}})$ are equivalent.

\strut

In [8], we discussed about it more in detail. In fact, under certain
condition, the converses also hold true. However, without the condition, we
cannot guarantee that the converses hold true.

\strut

\strut \strut \strut \strut

\strut

\section{Measure Theory on Graphs}

\strut

\strut

In this chapter, we will consider Measure theory on the given finite
directed graph $G,$ with respect to the reduced diagram graph measure space $%
(G^{\symbol{94}},$ $\mu _{G^{\symbol{94}}}).$ Of course, we can consider
Measure theory with respect to the graph measure space $(G,$ $\mu _{G}).$
But Calculus on $(G,$ $\mu _{G})$ can be regarded as the restricted Calculus
on $(G^{\symbol{94}},$ $\mu _{G^{\symbol{94}}}).$ So, we will concentrate on
observing Measure theory on $(G^{\symbol{94}},$ $\mu _{G^{\symbol{94}}}).$
\strut

\strut

\strut

\strut

\subsection{$G^{\symbol{94}}$-Measurable Functions}

\strut

\strut

Let $G$ be the given finite directed graph and $G^{\symbol{94}},$ the
shadowed graph and let $\Bbb{F}^{+}(G^{\symbol{94}})$ be the free
semigroupoid of $G^{\symbol{94}},$ and $(G^{\symbol{94}},$ $\mu _{G^{\symbol{%
94}}}),$ the $G^{\symbol{94}}$-measure space. All simple functions $g$ are
defined by

\strut

(3.1) $\ \ \ \ \ \ \ \ \ \ \ \ g=\sum_{n=1}^{N}a_{n}1_{S_{n}},$ \ for \ $%
a_{j}\in \Bbb{R},$

\strut

where $S_{1},$ ..., $S_{N}$ are subsets of $D_{r}(G^{\symbol{94}})$ and where

\strut

\begin{center}
$1_{S_{j}}(w)$ $=$ $\left\{ 
\begin{array}{ll}
1 & \text{if }w\in S_{j} \\ 
0 & \text{otherwise,}
\end{array}
\right. $
\end{center}

\strut

which are called the characteristic functions of $S_{j}$, for all $j$ $=$ $%
1, $ ..., $N.$ All characteristic functions are $G^{\symbol{94}}$%
-measurable. And all $G^{\symbol{94}}$-measurable functions are approximated
by simple functions.

\strut

\begin{example}
Let $w$ $\in $ $D_{r}(G^{\symbol{94}}).$ Then this element $w$ acts as one
of the $G^{\symbol{94}}$-measurable function from $D_{r}(G^{\symbol{94}})$
to $\Bbb{R}.$ Note that $w$ acts on $D_{r}(G^{\symbol{94}}),$ as the left
multiplication or the right multiplication. So, we have two different
functions $g_{l}^{w}$ and $g_{r}^{w}$ on $D_{r}(G^{\symbol{94}})$ such that

\strut

\ \ \ \ \ \ \ \ \ \ \ \ \ \ \ \ \ $\ g_{l}^{w}=1_{\delta ^{r}(S_{l}^{w})}$ \
\ and \ \ $g_{r}^{w}=1_{\delta ^{r}(S_{r}^{w})},$

\strut

for all $w^{\prime }\in D_{r}(G^{\symbol{94}}),$ where

\strut

\ $\ \ \ \ \ \ \ \ \ \ \ \ \ \ \ S_{l}^{w}=\{w^{\prime \prime }\in \Bbb{F}%
^{+}(G^{\symbol{94}}):ww^{\prime \prime }\in \Bbb{F}^{+}(G^{\symbol{94}})\}$

and

$\ \ \ \ \ \ \ \ \ \ \ \ \ \ \ S_{r}^{w}=\{w^{\prime \prime \prime }\in \Bbb{%
F}^{+}(G^{\symbol{94}}):w^{\prime \prime \prime }w\in \Bbb{F}^{+}(G^{\symbol{%
94}})\}.$

\strut \strut

Therefore, the element $w$ act as a $G^{\symbol{94}}$-measurable function $%
g_{w}$ on $D_{r}(G^{\symbol{94}})$ defined by

\strut

\ \ \ \ \ \ \ \ \ \ \ \ \ \ \ \ \ \ \ \ \ \ \ \ $\ \ \ \ g_{w}=1_{\delta
^{r}(S_{l}^{w})\cup \delta ^{r}(S_{r}^{w})}.$
\end{example}

\strut \strut \strut

\strut

\subsection{Integration on Graphs}

\strut

\strut

In this section, we will define the integrals of the given $G^{\symbol{94}}$%
-measurable functions, with respect to the $G^{\symbol{94}}$-measure $\mu
_{G^{\symbol{94}}}.$ First, let $g$ be a simple function given in (3.1).
Then the integral $I_{G}(g)$ of $g$ with respect to $\mu _{G^{\symbol{94}}}$
is defined by

\strut

(3.2) $\ \ \ \ \ \ \ \ \ I_{G}(g)\overset{denote}{=}\int_{G^{\symbol{94}%
}}g\,\,d\mu _{G^{\symbol{94}}}\overset{def}{=}\sum_{n=1}^{N}a_{n}\,\,\mu
_{G^{\symbol{94}}}(S_{n}).$

\strut

Since $\mu _{G^{\symbol{94}}}=d\cup \Delta ^{\symbol{94}},$ the definition
(3.2) can be rewritten as

\strut

\begin{center}
$
\begin{array}{ll}
I_{G}(g) & =\sum_{n=1}^{N}a_{n}\mu _{G^{\symbol{94}}}(S_{n})=%
\sum_{n=1}^{N}a_{n}\left( d(S_{n,V})+\Delta ^{\symbol{94}}(S_{n,FP})\right)
\\ 
&  \\ 
& =\sum_{n=1}^{N}a_{n}\left( \underset{v\in S_{n,V}}{\sum }\frac{\deg (v)}{%
\left| V(G)\right| }+\underset{\delta _{w}^{r}\in S_{n,FP}}{\sum }W(\delta
_{w}^{r})f(\delta _{w}^{r})\right) .
\end{array}
$
\end{center}

\strut \strut

\strut

\begin{proposition}
Let $g_{1}$ $=\sum_{j=1}^{n}1_{S_{j}}$ and $g_{2}=\sum_{i=1}^{m}1_{T_{i}}$
be simple functions, where $S_{j}$'s and $T_{i}$'s are subsets of $D_{r}(G^{%
\symbol{94}}).$ Suppose that $S_{j}$'s are mutually disjoint and also $T_{i}$%
's are mutually disjoint. If $\cup _{j=1}^{n}$ $S_{j}$ $=$ $\cup _{i=1}^{m}$ 
$T_{i}$ in $D_{r}(G^{\symbol{94}}),$ then $I_{G}(g_{1})$ $=$ $I_{G}(g_{2}).$
\end{proposition}

$\strut $

\begin{proof}
Observe that

\strut

(3.3)$\ \ \ \ \ I_{G}(g_{1})=\sum_{j=1}^{n}\left( \underset{v\in S_{j,V}}{%
\sum }\frac{\deg (v)}{\left| V(G)\right| }+\underset{w\in S_{j,FP}}{\sum }%
W(w)f(w)\right) $

and

(3.4)$\ \ \ \ \ I_{G}(g_{2})=\sum_{i=1}^{m}\left( \underset{v\in T_{i,V}}{%
\sum }\frac{\deg (v)}{\left| V(G)\right| }+\underset{w\in T_{i,FP}}{\sum }%
W(w)f(w)\right) .$

\strut

By the assumption, $\cup _{j=1}^{n}S_{j}$ and $\cup _{i=1}^{m}T_{i}$
coincide in $D_{r}(G^{\symbol{94}}).$ Therefore, we have that

\strut

\strut (3.5) \ \ $\ \ \ \ \ \ \ \cup _{j=1}^{n}S_{j}=\left( \cup
_{j=1}^{n}S_{j}\right) _{V}+\left( \cup _{j=1}^{n}S_{j}\right) _{FP}$

\strut

\ \ \ $\ \ \ \ \ \ \ \ \ \ \ \ \ \ \ \ \ \ \ \ \ \ \ \ \ \ \ =\left( \cup
_{i=1}^{m}T_{i}\right) _{V}+\left( \cup _{i=1}^{m}T_{i}\right) _{FP}$

\strut

\ \ \ $\ \ \ \ \ \ \ \ \ \ \ \ \ \ \ \ \ \ \ \ \ \ \ \ \ \ \ =\cup
_{i=1}^{m}T_{i}.$

\strut

Since $\{S_{j}\}_{j=1}^{n}$ and $\{T_{i}\}_{i=1}^{m}$ are disjoint families
of $P\left( D_{r}(G^{\symbol{94}})\right) ,$ the relation (3.5) explains

\strut

(3.6)\strut

\ \ \ \ \ \ \ \ \ \ \ \ \ $\ \ \ \ \ \ \ 
\begin{array}{ll}
\left( \cup _{j=1}^{n}S_{j}\right) _{V} & \text{\thinspace }=\cup
_{j=1}^{n}S_{j,V} \\ 
&  \\ 
& 
\begin{array}{l}
=\cup _{i=1}^{m}T_{i,V} \\ 
\\ 
=\left( \cup _{i=1}^{m}T_{i}\right) _{V}.
\end{array}
\end{array}
$

$\strut $

Similarly, we have that

\strut

(3.7)\ \ \ \ $\left( \cup _{j=1}^{n}S_{j}\right) _{FP}=\cup
_{j=1}^{n}S_{j,FP}=\cup _{i=1}^{m}T_{i,FP}=\left( \cup
_{i=1}^{m}T_{i}\right) _{FP}.$

\strut

Now, notice that

\strut

$\ \ \ \ \ \ \ \ \ \ \ \ \ \ \ \ \ \ \ \sum_{j=1}^{n}\underset{v\in S_{j,V}}{%
\sum }\frac{\deg (v)}{\left| V(G)\right| }=\underset{v\in \cup
_{j=1}^{n}S_{j,V}}{\sum }\frac{\deg (v)}{\left| V(G)\right| }$

and

$\ \ \ \ \ \ \ \ \ \ \ \ \ \ \ \ \ \ \ \sum_{i=1}^{m}\underset{v\in T_{i,V}}{%
\sum }\frac{\deg (v)}{\left| V(G)\right| }=\underset{v\in \cup
_{i=1}^{m}T_{i,V}}{\sum }\frac{\deg (v)}{\left| V(G)\right| },$

\strut

by the mutually disjointness on $\{S_{j,V}\}_{j=1}^{n}$ and $%
\{T_{i,V}\}_{i=1}^{m}.$ Similarly, by the mutually disjointness of $%
\{S_{j,FP}\}_{j=1}^{n}$ and $\{T_{i,FP}\}_{i=1}^{m}$ and by (3.7),

\strut

$\ \ \ \ \ \ \ \ \ \ \ \ \ \ \ \sum_{j=1}^{n}\underset{w\in S_{j,FP}}{\sum }%
L(\delta _{w}^{r})=\underset{w\in \cup _{j=1}^{n}S_{j,FP}}{\sum }L(w)$

and

$\ \ \ \ \ \ \ \ \ \ \ \ \ \ \ \sum_{i=1}^{m}\underset{w\in T_{i,FP}}{\sum }%
L(\delta _{w}^{r})=\underset{w\in \cup _{i=1}^{m}T_{i,FP}}{\sum }L(w).$

\strut \strut

Therefore, we can conclude that $I_{G}(g_{1})=I_{G}(g_{2}),$ by (3.5).
\end{proof}

\strut

It is easy to check that, if $g_{j}$ $=$ $\sum_{k=1}^{n}a_{j,k}1_{S_{j,k}},$
for $j$ $=$ $1,$ $2,$ then we can have

\strut

(3.8) $\ \ \ \ \ I_{G}\left( g_{1}+g_{2}\right) =I_{G}(g_{1})+I_{G}(g_{2})$

(3.9) \ \ \ $\ \ I_{G}(c\cdot g_{j})=c\cdot I_{G}(g_{j}),$ for all $j=1,2.$

\strut

Now, let $1_{S}$ and $1_{T}$ be the characteristic functions, where $S$ $%
\neq $ $T$ in $P\left( D_{r}(G^{\symbol{94}})\right) .$ Then we have that $%
1_{S}$ $\cdot $ $1_{T}$ $=$ $1_{S\cap T}.$

\strut \strut

\begin{proposition}
Let $g_{j}$ be given as above, for $j$ $=$ $1,$ $2.$ Then

\strut

$\ I_{G}(g_{1}g_{2})=\sum_{k,i=1}^{n}\left( a_{1,k}a_{2,i}\right) $

\strut

$\ \ \ \ \ \ \ \ \ \ \ \ \ \ \ \ \ \ \ \ \ \cdot \left( \underset{v\in
S_{(1,k),V}\cap S_{(2,i),V}}{\sum }\frac{\deg (v)}{\left| V(G)\right| }+%
\underset{w\in S_{(1,k),FP}\cap S_{(2,i),FP}}{\sum }L^{\symbol{94}%
}(w)\right) .$
\end{proposition}

\strut

\begin{proof}
Observe that

\strut

$\ \ \ g_{1}g_{2}=\left( \sum_{k=1}^{n}a_{1,k}1_{S_{1,k}}\right) \left(
\sum_{i=1}^{n}a_{2,i}1_{S_{2,i}}\right) $

\strut

$\ \ \ \ \ \ \ \ \ \ \ =\sum_{k,i=1}^{n}a_{1,k}a_{2,i}\left(
1_{S_{1,k}}\cdot 1_{S_{2,i}}\right)
=\sum_{k,i=1}^{n}a_{1,k}a_{2,i}1_{S_{1,k}\cap S_{2,i}}$

\strut

So, we have that

\strut

(3.10) $\ \ \ \ \ \ \ I_{G}(g_{1}g_{2})=\sum_{k,i=1}^{n}a_{1,k}a_{2,i}\mu
_{G^{\symbol{94}}}\left( S_{1,k}\cap S_{2,i}\right) .$

\strut

Consider

\strut

$\ \ \mu _{G^{\symbol{94}}}\left( S_{1,k}\cap S_{2,i}\right) =d\left(
(S_{1,k}\cap S_{2,i})_{V}\right) +\Delta ^{\symbol{94}}\left( (S_{1,k}\cap
S_{2,i})_{FP}\right) $

\strut

$\ \ \ \ \ \ \ \ \ \ \ =d\left( S_{(1,k),V}\cap S_{(2,i),V}\right) +\Delta ^{%
\symbol{94}}\left( S_{(1,k),FP}\cap S_{(2,k),FP}\right) $

\strut

$\ \ \ \ \ \ \ \ \ \ \ =\underset{v\in S_{(1,k),V}\cap S_{(2,i),V}}{\sum }%
\frac{\deg (v)}{\left| V(G)\right| }+\underset{w\in S_{(1,k),FP}\cap
S_{(2,i),FP})}{\sum }L^{\symbol{94}}(w).$

$\strut $
\end{proof}

\strut

Suppose that $g_{1}$ and $g_{2}$ are given as above and assume that the
families $\{S_{1,k}\}_{k=1}^{n}$ and $\{S_{2,i}\}_{i=1}^{n}$ are disjoint
family in $D_{r}(G^{\symbol{94}})$. i.e.,

$\strut $

\begin{center}
$\left( \cup _{k=1}^{n}S_{1,k}\right) \cap \left( \cup
_{i=1}^{n}S_{2,i}\right) =\varnothing .$
\end{center}

\strut

Then $I(g_{1}g_{2})=0.$

\strut

From now, similar to the classical measure theory, we can define the
integration of $G^{\symbol{94}}$-measurable functions.

\strut

\begin{definition}
Suppose that $f$ is a positive $G^{\symbol{94}}$-measurable function.. i.e., 
$f$ $=$ $f^{+}.$ Then we define the $G$-integral $I_{G}(f)$ of $f$ with
respect to $\mu _{G^{\symbol{94}}}$ by

\strut

$\ \ \ \ \ \ \ \ \ \ \ I_{G}(f)\overset{denote}{=}\int_{G^{\symbol{94}%
}}f\,\,d\mu _{G^{\symbol{94}}}=\underset{g\leq f\text{, }g\text{ is simple}}{%
\sup }I_{G}(g).$

\strut

Suppose that $f$ is $G^{\symbol{94}}$-measurable. Then, since $f$ $=$ $f^{+}$
$-$ $f^{-},$ we can define the $G$-integral of $f$ by

\strut

$\ \ \ \ \ \ \ \ \ \ \ \ \ \ \ \ \ \ \ \ \
I_{G}(f)=I_{G}(f^{+})-I_{G}(f^{-}).$
\end{definition}

\strut

Therefore, the positive $G^{\symbol{94}}$-measurable function $\left|
f\right| $ has its graph integral,

\strut

\begin{center}
$I_{G}\left( \left| f\right| \right) =I_{G}(f^{+})+I_{G}(f^{-}).$
\end{center}

\strut

If the $G^{\symbol{94}}$-measurable function $f$ satisfies $I_{G}(\left|
f\right| )$ $<$ $\infty ,$ then we say that this map $f$ is bounded.

\strut

Let's observe more $G^{\symbol{94}}$-measurable functions. If $g$ is a $G^{%
\symbol{94}}$-measurable function on $D_{r}(G^{\symbol{94}})$\strut , then
the support of $g$ is denoted by $D_{r}(G^{\symbol{94}}$ $:$ $g)$\strut .
First, consider the polynomial $g_{1}(x)$ $=$ $g_{x}$ on $D_{r}(G^{\symbol{94%
}}),$ where $g_{x}$ $=$ $1_{\delta ^{r}(S_{l}^{x})\cup \delta
^{r}(S_{r}^{x})}.$ This polynomial $g_{1}$ maps all $w$ in $\Bbb{F}^{+}(G^{%
\symbol{94}})$ to $g$\strut $_{w}$. Then, we can have the $G^{\symbol{94}}$%
-measurable function $g_{1}$ on $D_{r}(G^{\symbol{94}})$ defined by $%
g_{1}(w) $ $=$ $1_{\delta ^{r}(S_{l}^{w})\cup \delta ^{r}(S_{r}^{w})},$ for
all $w$ $\in $ $D_{r}(G).$

\strut

(3.11) $\ \ \ \ \ \ \ \ \ \ \ I_{G}(g_{1})=\underset{w\in D_{r}(G^{\symbol{94%
}})}{\sum }\mu _{G^{\symbol{94}}}\left( \delta ^{r}(S_{l}^{w})\cup \delta
^{r}(S_{r}^{w}))\right) ,$

\strut

by the previous proposition. Notice that the support $D_{r}(G^{\symbol{94}}$ 
$:$ $g_{1})$ of $g_{1}$ is identified with $D_{r}(G^{\symbol{94}}).$

\strut

Now, let $g_{2}(x)=g_{\delta _{x^{2}}^{r}}$ on $D_{r}(G^{\symbol{94}}).$
This polynomial $g_{2}$ has its support

\strut

\begin{center}
$D_{r}(G^{\symbol{94}}:g_{2})=\{w\in D_{r}(G^{\symbol{94}}):w\in V(G)$ or $w$
is a loop finite path$\}.$
\end{center}

\strut

The above support $D_{r}(G^{\symbol{94}}:g_{2})$ is determined by the fact
that $g_{2}(v)$ $=$ $g_{v^{2}}$ $=$ $g_{v},$ for all $v$ $\in $ $V(G),$ and $%
g_{2}(l)$ $=$ $g_{\delta _{l^{2}}^{r}}$ $=$ $g_{l},$ for all loop finite
paths $l$ in $D_{r}(G).$ Notice that the loop finite path in $D_{r}(G^{%
\symbol{94}})$ is the loop which is neither of the form $w^{-1}$ $w$ nor of
the form $w$ $w^{-1},$ by the reducing process (2.1). Assume now that $w$ $=$
$v_{1}$ $w$ $v_{2}$ is a non-loop finite path in $D_{r}(G^{\symbol{94}})$
with $v_{1}$ $\neq $ $v_{2}$ in $V(G^{\symbol{94}}).$ Then $w^{2}$ $=$ $%
(v_{1}$ $w$ $v_{2})$ $(v_{1}$ $w$ $v_{2})$ $=$ $\emptyset ,$ in $D_{r}(G^{%
\symbol{94}})$\strut $.$ Therefore, we can get that

\strut

(3.12) \ $I_{G}(g_{2})=\int_{G^{\symbol{94}}}g_{2}\,d\mu _{G^{\symbol{94}}}=%
\underset{w\in D_{r}(G^{\symbol{94}}:g_{2})}{\sum }\mu _{G^{\symbol{94}%
}}\left( \delta ^{r}(S_{l}^{w})\cup \delta ^{r}(S_{r}^{w})\right) .$

\strut

In general, we can get that;

\strut

\begin{proposition}
Let $g_{n}(x)=g_{x^{n}}$ be the monomial on $D_{r}(G^{\symbol{94}}).$ Then $%
I_{G}\left( g_{n}\right) $ $=$ $I_{G}\left( g_{2}\right) ,$ for all $n$ $\in 
$ $\Bbb{N}$ $\setminus $ $\{1\},$ where $I_{G}(g_{2})$ is given in (3.12).
\end{proposition}

\strut

\begin{proof}
It suffices to show that the support $D_{r}(G^{\symbol{94}}:g_{n})$ of $%
g_{n} $ and the support $D_{r}(G:g_{2})$ of $g_{2}$ coincide, for all $n$ $%
\in $ $\Bbb{N}$ $\setminus $ $\{1\}.$ It is easy to check that if $w$ is in $%
\Bbb{F}^{+}(G^{\symbol{94}}),$ then $w^{n}$ exists in $\Bbb{F}^{+}(G^{%
\symbol{94}})$ if and only if either $w$ is a vertex or $w$ is a loop finite
path. So, the support of $g_{n}$ is

\strut

$\ \ \ \ \ \ \ \ \ \ \ \ \ \ \ \ \ \ \ \ \ \ \ D_{r}(G^{\symbol{94}%
}:g_{n})=V(G)\cup loop_{r}(G),$

\strut

where $loop_{r}(G)\overset{def}{=}\{l\in D_{r}(G):l$ is a loop finite path$%
\}.$ Therefore, the support $D_{r}(G^{\symbol{94}}$ $:$ $g_{n})$ is same as $%
D_{r}(G^{\symbol{94}}$ $:$ $g_{2})$. Therefore,

\strut

$\ \ \ \ \ \ \ \ \ \ \ \ \ 
\begin{array}{ll}
I_{G}\left( g_{n}\right) & =\underset{w\in D_{r}(G^{\symbol{94}}:g_{2})}{%
\sum }I_{G}(g_{w}) \\ 
&  \\ 
& =\underset{w\in D_{r}(G^{\symbol{94}}:g_{2})}{\sum }\mu _{G^{\symbol{94}%
}}\left( \delta ^{r}(S_{l}^{w})\cup \delta ^{r}(S_{r}^{w})\right) .
\end{array}
$

$\strut \strut $
\end{proof}

\strut

\begin{theorem}
Let $g_{p}$ $=$ $\sum_{n=0}^{N}a_{n}g_{n}$ be a polynomial with $g_{0}$ $%
\equiv $ $1,$ for $a_{0},$ ..., $a_{N}$ $\in $ $\Bbb{R}.$ Then

\strut

\strut (3.13)

$\ \ 
\begin{array}{ll}
I_{G}\left( g_{p}\right) = & a_{0}\mu _{G^{\symbol{94}}}\left( D_{r}(G^{%
\symbol{94}})\right) +a_{1}\left( \underset{w\in D_{r}(G^{\symbol{94}})}{%
\sum }\mu _{G^{\symbol{94}}}\left( \delta ^{r}(S_{l}^{w})\cup \delta
^{r}(S_{r}^{w})\right) \right) \\ 
&  \\ 
& +\sum_{k=2}^{N}a_{k}\left( \underset{w\in V(G)\cup loop_{r}(G)}{\sum }\mu
_{G^{\symbol{94}}}\left( \delta ^{r}(S_{l}^{w})\cup \delta
^{r}(S_{r}^{w})\right) \right) ,
\end{array}
$

\strut

where $loop_{r}(G)=\{l\in D_{r}(G^{\symbol{94}}):l$ is a loop finite path$%
\}. $
\end{theorem}

\strut

\begin{proof}
Let $g_{p}$ be the given polynomial on $\Bbb{F}^{+}(G^{\symbol{94}}).$ Then
the integral $I_{G}\left( g_{p}\right) $ of $g_{p}$ is determined by;

\strut

$\ \ \ I_{G}\left( g_{p}\right) =\sum_{n=0}^{N}a_{n}I_{G}\left( g_{n}\right) 
$

\strut

$\ \ \ \ \ \ \ =a_{0}I_{G}(1)+a_{1}I_{G}\left( g_{1}\right)
+\sum_{k=2}^{N}a_{k}I_{G}\left( g_{k}\right) $

\strut

$\ \ \ \ \ \ \ =a_{0}\mu _{G^{\symbol{94}}}\left( D_{r}(G^{\symbol{94}%
})\right) +a_{1}\left( \underset{w\in D_{r}(G^{\symbol{94}})}{\sum }\mu _{G^{%
\symbol{94}}}(\delta ^{r}(S_{l}^{w})\cup \delta ^{r}(S_{r}^{w}))\right) $

\strut

$\ \ \ \ \ \ \ \ \ \ \ \ \ \ \ \ \ \ \ \ \ \ \ \ \ \ \ \ \ \ \ \ \ \ \ \ \ \
\ +\sum_{k=2}^{N}a_{k}I_{G}\left( g_{2}\right) $

\strut

since the constant function $1$ has its support, $D_{r}(G^{\symbol{94}}),$
and hence, since $I_{G}\left( g_{n}\right) $ $=$ $I_{G}\left( g_{2}\right) ,$
for all $n$ $\geq $ $2,$ we can get that

\strut

$\ \ \ \ \ \ \ =a_{0}\mu _{G^{\symbol{94}}}\left( D_{r}(G)\right)
+a_{1}\left( \underset{w\in D_{r}(G^{\symbol{94}})}{\sum }\mu _{G^{\symbol{94%
}}}(\delta ^{r}(S_{l}^{w})\cup \delta ^{r}(S_{r}^{w}))\right) $

\strut

$\ \ \ \ \ \ \ \ \ \ \ \ \ \ \ \ \ \ \ \ \ \ \ \ \ \ \ \ \ \ \
+\sum_{k=2}^{N}a_{k}\left( \underset{w\in V(G)\cup loop(G^{\symbol{94}})}{%
\sum }\mu _{G^{\symbol{94}}}\left( \delta ^{r}(S_{l}^{w})\cup \delta
^{r}(S_{r}^{w})\right) \right) ,$

\strut

because the support $D_{r}(G^{\symbol{94}}:g_{2})$ of $g_{2}$ is the union
of $V(G^{\symbol{94}})$ and $loop_{r}(G^{\symbol{94}}),$ where $loop_{r}(G^{%
\symbol{94}})$ is the collection of all loop finite paths in $D_{r}(G^{%
\symbol{94}}).$
\end{proof}

\strut \strut \strut

\begin{corollary}
Let $w$ $=v_{1}wv_{2}$ be a finite path in $D_{r}(G^{\symbol{94}})$ with $%
v_{1}$ $\neq $ $v_{2}$ in $V(G)$, and let $g$ $=$ $\sum_{n=0}^{N}$ $%
a_{n}g_{w}^{n},$ with $g_{w}^{0}$ $\overset{def}{=}$ $1,$ where $a_{0},$
..., $a_{N}$ $\in $ $\Bbb{R},$ and $g_{w}$ is $1_{\delta ^{r}(S_{l}^{w})\cup
\delta ^{r}(S_{r}^{w})}.$ Then

\strut

$\ \ \ \ \ \ \ I_{G}\left( g\right) =a_{0}\mu _{G^{\symbol{94}}}\left(
D_{r}(G^{\symbol{94}})\right) +a_{1}\mu _{G^{\symbol{94}}}\left( \delta
^{r}(S_{l}^{w})\cup \delta ^{r}(S_{r}^{w})\right) .$
\end{corollary}

\strut

\begin{proof}
Since $w$ is a non-loop finite path, $w^{k}$ $=$ $\emptyset ,$ for all $k$ $%
\in $ $\Bbb{N}$ $\setminus $ $\{1\}.$ Therefore, $S_{l}^{w^{k}}$ $\cup $ $%
S_{r}^{w^{k}}$ $=$ $\varnothing ,$ for all $k$ $=$ $2,$ $3,$ ..., $N.$ This
shows that $I_{G}\left( g_{w}^{k}\right) $ $=$ $0,$ for all $k$ $=$ $2,$ $3,$
..., $N.$ So, $I_{G}\left( g\right) $ $=$ $a_{0}$ $I_{G}\left( 1\right) $ $+$
$a_{1}$ $I_{G}\left( g_{w}\right) .$
\end{proof}

\strut

\begin{corollary}
Let $w$ $=$ $vwv$ be a loop finite path in $D_{r}(G^{\symbol{94}}),$ with $v$
$\in $ $V(G^{\symbol{94}}),$ and let $g$ $=$ $\sum_{n=0}^{N}$ $%
a_{n}g_{w}^{n} $, with $g_{w}^{0}$ $=$ $1,$ where $a_{0},$ ..., $a_{n}$ $\in 
$ $\Bbb{R}.$ Then

\strut

$\ \ \ \ \ \ I_{G}\left( g\right) =a_{0}\mu _{G^{\symbol{94}}}\left(
D_{r}(G)\right) +\sum_{k=1}^{N}a_{k}\mu _{G^{\symbol{94}}}\left( \delta
^{r}(S_{l}^{w})\cup \delta ^{r}(S_{r}^{w})\right) .$

$\square $
\end{corollary}

\strut \strut \strut \strut

Now, we will consider the map $g_{-1}$ defined by $g_{-1}(x)$ $=$ $%
g_{x^{-1}} $ $=$ $g_{\delta _{x^{-1}}^{r}}.$ On $D_{r}(G^{\symbol{94}}),$
the map $g_{-1}$ is well-determined, because, for any $w$ $\in $ $\Bbb{F}%
^{+}(G^{\symbol{94}}),$ there always exists $w^{-1}.$ Also, the
admissibility of $w^{-1}$ is preserved by that of $w.$ This discussion says
that the support $D_{r}(G^{\symbol{94}}$ $:$ $g_{\delta _{w}^{r}})$ of $%
g_{\delta _{w}^{r}}$ and the support $D_{r}(G^{\symbol{94}}$ $:$ $g_{\delta
_{w^{-1}}^{r}})$ of $g_{\delta _{w^{-1}}^{r}}$ are same. i.e., If $w$ is in $%
D_{r}(G^{\symbol{94}}),$ then

\strut

(3.14) \ \ \ $D_{r}(G^{\symbol{94}}:g_{w})=D_{r}(G^{\symbol{94}%
}:g_{w^{-1}}), $ for all $w\in D_{r}(G^{\symbol{94}}).$

\strut

Futhermore, we can get that

\strut

(3.15)\ \ \ \ \ \ \ \ \ \ \ $D_{r}(G^{\symbol{94}}:g_{-1})=D_{r}(G^{\symbol{%
94}})=D_{r}(G^{\symbol{94}}:g_{1}),$

\strut

where $g_{1}(x)=g_{x},$ for all $x$ $\in $ $D_{r}(G^{\symbol{94}})$\strut ,
since $g_{-1}$ has its support $D_{r}(G^{\symbol{94}})$.

\strut

\begin{proposition}
Let $g_{-1}$ be given as above. Then $g_{-1}$ is $G^{\symbol{94}}$%
-measurable and $I_{G}\left( g_{-1}\right) $ $=$ $I_{G}\left( g_{1}\right) .$
$\square $
\end{proposition}

\strut \strut \strut

Similarly, we can conclude that;

\strut

\begin{proposition}
Define $g_{-n}(x)=g_{\delta _{x^{-n}}^{r}},$ for $n\in \Bbb{N}$, for $x$ $%
\in $ $D_{r}(G^{\symbol{94}}).$ Then $I_{G}\left( g_{-n}\right) =I_{G}\left(
g_{n}\right) $.
\end{proposition}

\strut

\begin{proof}
Observe that if $n$ $\geq $ $2,$ then

\strut

$\ \ D_{r}(G^{\symbol{94}}:g_{-n})=D_{r}(G^{\symbol{94}}:g_{-2})=V(G^{%
\symbol{94}})\cup loop_{r}(G^{\symbol{94}})=D_{r}(G^{\symbol{94}}:g_{2}).$

\strut

So, $I_{G}\left( g_{-n}\right) =I_{G}\left( g_{-2}\right) =I_{G}\left(
g_{2}\right) .$ By the previous proposition, $I_{G}\left( g_{-1}\right) $ $=$
$I_{G}\left( g_{1}\right) .$
\end{proof}

\strut

So, we can consider the trigonometric polynomials on $\Bbb{F}^{+}(G^{\symbol{%
94}}).$

\strut

\begin{theorem}
Let $g_{k}$ and $g_{-k}$ be given as above, for all $k$ $\in $ $\Bbb{N},$
and let $g$ $=$ $\sum_{n=-M}^{N}$ $a_{n}g_{n}$ be a trigonometric
polynomial, with $g_{0}$ $\equiv $ $1,$ where $N,$ $M$ $\in $ $\Bbb{N}.$ Then

\strut

(3.16) $\ I_{G}\left( g\right) =a_{0}\mu _{G^{\symbol{94}}}\left(
D_{r}(G)\right) +\left( a_{1}+a^{-1}\right) I_{G}(g_{1})$

\strut

$\ \ \ \ \ \ \ \ \ \ \ \ \ \ \ \ \ \ \ \ \ \ \ \ \ \ \ \ \ \ \ \
+\sum_{n=-M}^{-2}a_{n}I_{G}\left( g_{2}\right)
+\sum_{k=2}^{N}a_{k}I_{G}\left( g_{2}\right) .$
\end{theorem}

\strut

\begin{proof}
Let $g$ be the given trigonometric polynomial on $D_{r}(G^{\symbol{94}}).$
Then

\strut

$I_{G}\left( g\right) =I_{G}\left( \sum_{n=-N}^{N}a_{n}g_{n}\right)
=\sum_{n=-N}^{N}a_{n}I_{G}\left( g_{n}\right) $

\strut

\ $=\sum_{n=-M}^{-2}a_{n}I_{G}\left( g_{n}\right)
+a_{-1}I_{G}(g_{-1})+a_{0}I_{G}\left( g_{0}\right)
+a_{1}I_{G}(g_{1})+\sum_{k=2}^{N}a_{k}I_{G}\left( g_{k}\right) $

\strut

\ $=\sum_{n=-M}^{-2}a_{n}I_{G}\left( g_{n}\right) +a_{0}\mu _{G^{\symbol{94}%
}}\left( D_{r}(G)\right) ++\sum_{k=1}^{N}a_{k}I_{G}\left( g_{k}\right) $

\strut

by the previous proposition

\strut

$=a_{0}\mu _{G^{\symbol{94}}}\left( D_{r}(G)\right)
+(a_{-1}+a_{1})I_{G}(g_{1})+\sum_{n=-M}^{-1}a_{n}I_{G}\left( g_{2}\right)
+\sum_{k=1}^{N}a_{k}I_{G}\left( g_{2}\right) $

\strut

by the fact that $I_{G}\left( g_{n}\right) =I_{G}\left( g_{2}\right) ,$ for
all $n$ $\in $ $\Bbb{N}$ $\setminus $ $\{1\}.$
\end{proof}

\strut \strut \strut

\strut \strut

\strut \strut

\subsection{Examples}

\strut

\strut

In this section, we will consider certain finite directed graphs. Let $%
G_{\Lambda }$ be a tree with

\strut

\begin{center}
$V(G_{\Lambda })=\{v_{1},v_{2},v_{3}\}$ \ and \ $E(G_{\Lambda })=\{e_{1}$ $=$
$v_{1}e_{1}v_{2},$ $e_{2}=v_{1}e_{2}v_{3}\}.$
\end{center}

\strut

Let $G_{\Delta }$ be a circulant graph with

\strut

\begin{center}
$V(G_{\Delta })=\{v_{1},v_{2},v_{3}\}$ \ and \ $E(G_{\Delta })=\left\{ 
\begin{array}{c}
e_{1}=v_{1}e_{1}v_{2}, \\ 
e_{2}=v_{2}e_{2}v_{3}, \\ 
e_{3}=v_{3}e_{3}v_{1}
\end{array}
\right\} .$
\end{center}

\strut

Also, assume that the corresponding reduced diagram graph measures are
unweighted. i.e.,

\strut

\begin{center}
$W_{\Lambda }=1=W_{\Delta },$
\end{center}

\strut

where $W_{\Lambda }$ and $W_{\Delta }$ are the weighting map on $%
FP(G_{\Lambda }^{\symbol{94}})$ and $FP(G_{\Delta }^{\symbol{94}}),$
respectively.

\strut

\begin{example}
Consider $G_{\Lambda }.$ We have that

\strut

$\ \ 
\begin{array}{ll}
I_{G_{\Lambda }}\left( g_{v_{1}}\right) & =\mu _{G_{\Lambda }^{\symbol{94}%
}}\left( \{v_{1},e_{1},e_{2},e_{1}^{-1},e_{2}^{-1}\}\right)
=d(\{v_{1}\})+\Delta ^{\symbol{94}}\left(
\{e_{1},e_{2},e_{1}^{-1},e_{2}^{-1}\}\right) \\ 
& =d(v_{1})+\Delta ^{\symbol{94}}\left( e_{1}\right) +\Delta ^{\symbol{94}%
}(e_{2})+\Delta ^{\symbol{94}}\left( e_{1}^{-1}\right) +\Delta ^{\symbol{94}%
}\left( e_{2}^{-1}\right) \\ 
& =\frac{4}{3}+4=\frac{16}{3}.
\end{array}
$

\strut

$\ \ \ \ \ \ \ \ \ \ \ \ \ \ \ 
\begin{array}{ll}
I_{G_{\Lambda }}\left( g_{v_{2}}\right) & =\mu _{G_{\Lambda }^{\symbol{94}%
}}(\{v_{2},e_{1},e_{1}^{-1}\})=d(\{v_{2}\})+\Delta ^{\symbol{94}%
}(\{e_{1},e_{1}^{-1}\}) \\ 
& =\frac{2}{3}+2=\frac{8}{3}.
\end{array}
$

\strut

Similarly, $I_{G_{\Lambda }}\left( g_{v_{3}}\right) =\frac{8}{3}.$

\strut

$\ \ \ \ \ \ \ \ \ \ \ \ \ \ \ 
\begin{array}{ll}
I_{G_{\Lambda }}\left( g_{e_{1}}\right) & =\mu _{G_{\Lambda }^{\symbol{94}%
}}\left( \{v_{1},v_{2},e_{2}^{-1}\}\right) =d\left( \{v_{1},v_{2}\}\right)
+\Delta ^{\symbol{94}}(e_{2}^{-1}) \\ 
& =d(v_{1})+d(v_{2})+1=\frac{\deg (v_{1})}{\left| V(G)\right| }+\frac{\deg
(v_{2})}{\left| V(G)\right| }+1 \\ 
& =\frac{4}{3}+\frac{2}{3}+1=3=I_{G_{\Lambda }}\left( g_{e_{1}^{-1}}\right) .
\end{array}
$

\strut

Similarly, $I_{G_{\Lambda }}\left( g_{e_{2}}\right) =3=I_{G_{\Lambda
}}\left( g_{e_{2}^{-1}}\right) .$ We have that

\strut

$\ \ \ \ \ \ \ \ \ \ \ \ \ \ \ \ \ D_{r}(G_{\Lambda
})=\{v_{1},v_{2},v_{3}\}\cup \{e_{1},e_{2},e_{1}^{-1},e_{2}^{-1}\}.$

\strut

Thus

\strut

$\ \ \ \ 
\begin{array}{ll}
I_{G_{\Lambda }}\left( g_{1}\right) & =\underset{w\in D_{r}(G_{\Lambda }^{%
\symbol{94}})}{\sum }I_{G_{\Lambda }}\left( g_{w}\right) \\ 
& =\sum_{j=1}^{3}I_{G_{\Lambda }}\left( g_{v_{j}}\right)
+\sum_{k=1}^{2}I_{G_{\Lambda }}\left( g_{e_{k}}\right)
+\sum_{i=1}^{2}I_{G_{\Lambda }}\left( g_{e_{i}^{-1}}\right) \\ 
& =\left( \frac{16}{3}+\frac{8}{3}+\frac{8}{3}\right) +\left( 2+2\right)
+\left( 2+2\right) =\allowbreak \frac{56}{3}.
\end{array}
$

$\strut \strut $

We also have that

\strut

$\ \ \ \ \ \ \ \ \ I_{G_{\Lambda }}\left( g_{n}\right)
=\sum_{j=1}^{3}I_{G_{\Lambda }}\left( g_{v_{j}}\right) =\frac{16}{3}+\frac{8%
}{3}+\frac{8}{3}=\allowbreak \frac{32}{3},$

\strut

for all $n$ $\in $ $\Bbb{N}$ $\setminus $ $\{1\},$ because $w^{k}$ $=$ $%
\emptyset ,$ for all $w$ $\in $ $D_{r}^{FP}(G_{\Lambda }^{\symbol{94}})$ and 
$k$ $\in $ $\Bbb{N}$ $\setminus $ $\{1\},$ and $D_{r}(G_{\Lambda }^{\symbol{%
94}}$ $:$ $g_{k})$ $=$ $V(G_{\Lambda }^{\symbol{94}})$ $\cup $ $%
loop_{r}(G_{\Lambda }^{\symbol{94}})$ $=$ $V(G_{\Lambda }^{\symbol{94}}).$
Therefore,

\strut

$\ \ \ \ \ \ \ \ \ \ \ \ \ \ \ \ \ \ \ \ \ \ \ \ \ \ \ \ I_{G_{\Lambda
}}\left( g_{p}\right) =\frac{56}{3}+\frac{32(N-1)}{3}.$

\strut

where $g_{p}=\sum_{n=1}^{N}g_{n}$, with $g_{0}$ $=$ $0,$ for all $N$ $\in $ $%
\Bbb{N}.$
\end{example}

\strut

\begin{example}
Consider the graph $G_{\Delta }.$ Then we have;

\strut

$\ \ \ 
\begin{array}{ll}
I_{G_{\Delta }}\left( g_{v_{1}}\right) & =\mu _{G_{\Delta }^{\symbol{94}%
}}\left( \{v_{1},e_{1}^{\pm 1},e_{3}^{\pm 1}\}\right) =d(\{v_{1}\})+\Delta ^{%
\symbol{94}}\left( \{e_{1}^{\pm 1},e_{3}^{\pm 1}\}\right) \\ 
& =\frac{\deg (v)}{\left| V(G)\right| }+2L\left( e_{1}\right) +2L(e_{3})=%
\frac{4}{3}+2+2 \\ 
& =\frac{16}{3}.
\end{array}
$

\strut

Similarly, $I_{G_{\Delta }}\left( g_{v_{2}}\right) =\frac{16}{3}%
=I_{G_{\Delta }}\left( g_{v_{3}}\right) .$ By the existence of the circulant
parts, we can have that

\strut

\ \ \ \ \ \ \ \ \ \ \ \ $I_{G_{\Delta }}\left( g_{e_{j}}\right) =\frac{17}{3}%
=I_{G_{\Delta }}\left( g_{e_{j}^{-1}}\right) ,$ for all $j$ $=$ $1,$ $2,$ $%
3. $

\strut

Let $j$ $=$ $1.$ Then $I_{G_{\Delta }}\left( g_{e_{1}}\right) =\mu
_{G_{\Delta }^{\symbol{94}}}\left( \delta ^{r}(S_{l}^{w})\cup \delta
^{r}(S_{r}^{w})\right) ,$ where

\strut

$\ \ \ \ \ \ \ \ \ \ \ \ \ \ \ \ S_{l}^{w}=\{v_{2},e_{2}\}$ \ and \ $%
S_{r}^{w}=\{v_{1},e_{1}^{-1},e_{3}\},$

so,

$\ \ \ \ \ \ \ \ \ \ \ \ \ \ \ \ \ \ \ \ S_{l}^{w}\cup
S_{r}^{w}=\{v_{1},v_{2},e_{1}^{-1},e_{2},e_{3}\}.$

\strut Thus

\strut

$\ \ \ I_{G_{\Delta }}\left( g_{e_{1}}\right) =\mu _{G_{\Delta }^{\symbol{94}%
}}\left( \{v_{1},v_{2},e_{1}^{-1},e_{2},e_{3}\}\right) $

\strut

$\ \ \ \ \ \ \ \ \ \ \ \ \ \ \ \ \ \ =d\left( \{v_{1},v_{2}\}\right) +\Delta
^{\symbol{94}}\left( \{e_{1}^{-1},e_{2},e_{3}\}\right) $

\strut

$\ \ \ \ \ \ \ \ \ \ \ \ \ \ \ \ \ \ =\frac{\deg (v_{1})}{\left| V(G)\right| 
}+\frac{\deg (v_{2})}{\left| V(G)\right| }+L\left( e_{1}^{-1}\right)
+L(e_{2})+L(e_{3})$

\strut

$\ \ \ \ \ \ \ \ \ \ \ \ \ \ \ \ \ \ =\frac{4}{3}+\frac{4}{3}%
+1+1+1=\allowbreak \frac{17}{3}.$

\strut

Now, consider the element $l=e_{1}e_{2}e_{3}$ in $D_{r}(G^{\symbol{94}}).$
Then

\strut

$\ \ \ \ \ \ \ \ \ \ \ \ \ \ \delta ^{r}(S_{l}^{l})=\{v_{1},l\}$ \ and \ $%
\delta ^{r}(S_{r}^{l})=\{v_{1},l\},$

so,

\ \ \ $\ \ \ \ \ \ \ \ \ \ \ \ \ \ \ \ \ \ \ \delta ^{r}(S_{l}^{l})\cup
\delta ^{r}(S_{r}^{l})=\{v_{1},l\}.$

\strut Therefore,

\strut

$\ \ \ \ \ \ \ \ \ \ \ \ \ \ I_{G_{\Delta }}\left( g_{l}\right) =\mu _{G^{%
\symbol{94}}}\left( \{v_{1},l\}\right) =\frac{\deg (v_{1})}{\left|
V(G)\right| }+L\left( l\right) $

\strut \strut

$\ \ \ \ \ \ \ \ \ \ \ \ \ \ \ \ \ =\frac{4}{3}+3=\allowbreak \frac{13}{3}.$

\strut \strut \strut \strut \strut
\end{example}

\strut \strut \strut \strut

\strut

\section{Subgraph Measure Theory}

\strut

\strut

\strut

In this chapter, we will consider a new measure on the finite directed
graph. Throughout this chapter, let $G$ be a finite directed graph and $G^{%
\symbol{94}},$ the corresponding shadowed graph and let $\Bbb{F}^{+}(G^{%
\symbol{94}})$ be the free semigroupoid of $G^{\symbol{94}}.$ Suppose $H$ is
a full subgraph of $G.$ The shadowed graph $H^{\symbol{94}}$ of $H$ is also
a subgraph of the shadowed graph $G^{\symbol{94}}$ of $G.$ So, canonically,
we can construct the $H$-measures $\mu _{H}$, as in Chapter 1, and the $H^{%
\symbol{94}}$-measure $\mu _{H^{\symbol{94}}},$ as in Chapter 2. We restrict
our interests to the $H^{\symbol{94}}$-measure $\mu _{H^{\symbol{94}}}.$

\strut

\begin{definition}
Let $(G^{\symbol{94}},$ $\mu _{G^{\symbol{94}}})$ be the reduced diagram
graph measure space of $G$ and let $H$ be a full subgraph of $G.$ Define the 
$H^{\symbol{94}}$-measure $\mu _{H^{\symbol{94}}}$ by $d_{H}$ $\cup $ $%
\Delta _{H}^{\symbol{94}},$ where

\strut

(4.1) $\ \ \ \ \ \ \ \ \ \ \ \ \ \ \ \ \ \ \ \ \ d_{H}(S)\overset{def}{=}%
\underset{v\in V(H)\cap S_{V}}{\sum }\frac{\deg (v)}{\left| V(H)\right| },$

\strut 

for all $S_{V}\subseteq V(G^{\symbol{94}})$, and

\strut 

(4.2) $\ \ \ \Delta _{H}^{\symbol{94}}(S_{FP})\overset{def}{=}\underset{w\in
FP(H)\cap S_{FP}}{\sum }L(w)=\underset{w\in FP(H)\cap S_{FP}}{\sum }W(\delta
_{w}^{r})f\left( \delta _{w}^{r}\right) ,$

\strut 

for all $S_{FP}\subseteq D_{r}^{FP}(G).$
\end{definition}

\strut

By definition, the subgraph measure $\mu _{H^{\symbol{94}}}$ is the
restricted measure $\mu _{G^{\symbol{94}}}\mid _{D_{r}(H^{\symbol{94}})}.$

\strut

\begin{proposition}
Let $g_{w}$ $=$ $1_{\delta ^{r}(S_{l}^{w})\cup \delta ^{r}(S_{r}^{w})}$ be
the $G^{\symbol{94}}$-measurable function defined in Section 3.2, for each $w
$ $\in $ $D_{r}(G^{\symbol{94}}).$ Then graph integral $I_{G:H}\left(
g_{w}\right) $ of $g_{w}$ with respect to the subgraph measure $\mu _{H^{%
\symbol{94}}}$ is determined by

\strut 

$\ \ \ \ \ \ \ 
\begin{array}{ll}
I_{G:H}\left( g_{w}\right)  & =\mu _{H^{\symbol{94}}}\left( \delta
^{r}(S_{l}^{w})\cup \delta ^{r}(S_{r}^{w})\right)  \\ 
&  \\ 
& =\mu _{G^{\symbol{94}}}\left( D_{r}(H^{\symbol{94}})\cap (\delta
^{r}(S_{l}^{w})\cup \delta ^{r}(S_{r}^{w}))\right) .
\end{array}
$

$\square $
\end{proposition}

\strut \strut \strut \strut \strut \strut \strut \strut

\strut

\strut

\section{Extended Graph Measures}

\strut

\strut

In this chapter, based on the graph measure construction in Chapter 2, we
will consider the extended graph measuring. Throughout this chapter, let $G$
be a finite directed graph and $G^{\symbol{94}},$ the shadowed graph and let 
$\Bbb{F}^{+}(G^{\symbol{94}})$ and $D_{r}(G^{\symbol{94}})$ be the
corresponding free semigroupoid of $G^{\symbol{94}}$ and the reduced diagram
set of $G^{\symbol{94}},$ respectively. As usual, let $\delta ^{r}$ $:$ $%
\Bbb{F}^{+}(G^{\symbol{94}})$ $\rightarrow $ $D_{r}(G^{\symbol{94}})$ be the
reduced diagram map which is a surjection. By definition, the power set $%
P\left( \Bbb{F}^{+}(G^{\symbol{94}})\right) $ is a $\sigma $-algebra of $%
\Bbb{F}^{+}(G^{\symbol{94}}).$ Indeed,

\strut

(i) \ \ $\varnothing ,$ $\Bbb{F}^{+}(G^{\symbol{94}})\in P\left( \Bbb{F}%
^{+}(G^{\symbol{94}})\right) $

(ii) \ If $S$ $\subseteq $ $\Bbb{F}^{+}(G^{\symbol{94}}),$ then $\Bbb{F}%
^{+}(G^{\symbol{94}})$ $\setminus $ $S$ $\subseteq $ $\Bbb{F}^{+}(G^{\symbol{%
94}})$

(iii) If $(S_{k})_{k=1}^{\infty }$ $\subseteq $ $P\left( \Bbb{F}^{+}(G^{%
\symbol{94}})\right) ,$ then $\cup _{k=1}^{\infty }S_{k}$ $\subseteq $ $\Bbb{%
F}^{+}(G^{\symbol{94}}).$

\strut

We want to define a suitable measure on $\Bbb{F}^{+}(G^{\symbol{94}}).$

\strut

\begin{definition}
Define a measure $\mu $ on $\Bbb{F}^{+}(G^{\symbol{94}})$ by $d$ $\cup $ $%
\omega ,$ where $d$ is the weighted degree measure on $V(G^{\symbol{94}})$
defined in Chapter 2. i.e.,

\strut 

$\ \ \ \ \ \ \ \ \ \ \ \ \ \ \ \ \ \ \ \ d\left( S_{V}\right) $ $=$ $%
\underset{v\in S_{V}}{\sum }\frac{\deg (v)}{\left| V(G^{\symbol{94}})\right| 
},$

\strut 

and $\omega $ is the measure on $FP(G^{\symbol{94}})$ defined by

\strut 

$\ \ \ \ \ \ \ \ \ \ \ \ \ \ \ \omega (S_{FP})=\underset{w\in S_{FP}}{\sum }%
W(\delta _{w}^{r})f(\delta _{w}^{r}),$

\strut 

where $W$ and $f$ are the weighting map and the reduced diagram length map
defined on $D_{r}(G^{\symbol{94}}),$ in Chapter 2. The measure $\mu $ is
called the (extended) graph measure on $\Bbb{F}^{+}(G^{\symbol{94}}).$
\end{definition}

\strut

It is easy to check that $\mu $ is a measure on $\Bbb{F}^{+}(G^{\symbol{94}%
}).$ Remark the difference between $\mu _{G^{\symbol{94}}}$ on $D_{r}(G^{%
\symbol{94}})$ and $\mu $ on $\Bbb{F}^{+}(G^{\symbol{94}}).$ A difference
between $\mu $ and $\mu _{G^{\symbol{94}}}$ is recognized by the following
example. Suppose $l$ is a loop-edge in $E(G^{\symbol{94}})$ and let $S_{1}$ $%
=$ $\{l\}$ and $S_{2}$ $=$ $\{l^{2}\}.$ Then the union $S$ $=$ $S_{1}$ $\cup 
$ $S_{2}$ is in $\Bbb{F}^{+}(G^{\symbol{94}}),$ but $S$ is not in $D_{r}(G^{%
\symbol{94}}).$ So, the measure $\mu _{G^{\symbol{94}}}(S)$ is undefined.
However, it is defined on $\delta ^{r}(S)$ $=$ $\{\delta _{l}^{r},$ $\delta
_{l^{2}}^{r}\}$ $=$ $\{l\}$. Moreover, we have that

\strut

\begin{center}
$
\begin{array}{ll}
\mu _{G^{\symbol{94}}}\left( \delta ^{r}(S)\right)  & =\mu _{G^{\symbol{94}%
}}\left( \delta ^{r}(S_{1})\cup \delta ^{r}(S_{2})\right)  \\ 
&  \\ 
& =\Delta ^{\symbol{94}}\left( \delta ^{r}(S_{1})\cup \delta
^{r}(S_{2})\right)  \\ 
&  \\ 
& =\Delta ^{\symbol{94}}\left( \{l\}\cup \{l\}\right) =\Delta ^{\symbol{94}%
}(\{l\})=W(l)
\end{array}
$
\end{center}

\strut and

\begin{center}
\strut $
\begin{array}{ll}
\mu \left( S\right)  & =\mu (S_{1}\cup S_{2})=\mu (\{l,l^{2}\})=\omega
\left( \{l,l^{2}\}\right)  \\ 
&  \\ 
& =\omega (\{l\})+\omega (\{l=\delta _{l^{2}}^{r}\})=2W(l).
\end{array}
$
\end{center}

\strut \strut

This difference shows us that the measure $\mu $ is not bounded, in general.
However, this measure $\mu $ is locally bounded in the sense that $\mu (S)$ $%
<$ $\infty ,$ for all finite subsets $S$ in $\Bbb{F}^{+}(G^{\symbol{94}}).$

\strut

\begin{proposition}
The extended graph measure $\mu $ on $\Bbb{F}^{+}(G^{\symbol{94}})$ is a
locally bounded positive measure. $\square $
\end{proposition}

\strut

\begin{definition}
The triple $\left( \Bbb{F}^{+}(G^{\symbol{94}}),\text{ }P\left( \Bbb{F}%
^{+}(G^{\symbol{94}})\right) ,\text{ }\mu \right) $ is called the extended
graph measure space (in short, $\mu $-measure space). For convenience, we
denote it by $(G^{\symbol{94}},$ $\mu ).$
\end{definition}

\strut

In the following section, we will show that the similar integration theory
holds for the extended graph measure like in Chapter 3.\strut

\strut \strut

\strut

\strut

\subsection{\strut Integration on Graphs}

\strut

\strut

In this section, we will define the integrals of the given $\mu $-measurable
functions, with respect to the extended graph measure $\mu .$ First, let

$\strut $

(5.1) $\ \ \ \ \ \ \ \ \ \ \ g=\sum_{n=1}^{N}a_{n}1_{S_{n}},$ for all $%
S_{n}\subseteq \Bbb{F}^{+}(G^{\symbol{94}}),$

\strut

where $a_{1},$ ..., $a_{N}$ $\in $ $\Bbb{R},$ be a simple function. Then the
integral $I_{\mu }(g)$ of $g$ with respect to $\mu $ is defined by

\strut

(5.2) $\ \ \ \ \ \ \ \ \ I_{G}(g)\overset{denote}{=}\int_{G^{\symbol{94}%
}}g\,\,d\mu \overset{def}{=}\sum_{n=1}^{N}a_{n}\,\,\mu (S_{n}).$

\strut

Since $\mu =d\cup \omega ,$ the definition (5.2) can be rewritten as

\strut

\begin{center}
$
\begin{array}{ll}
I_{G}(g) & =\sum_{n=1}^{N}a_{n}\mu _{G^{\symbol{94}}}(S_{n})=%
\sum_{n=1}^{N}a_{n}\left( d(S_{n,V})+\omega (S_{n,FP})\right) \\ 
&  \\ 
& =\sum_{n=1}^{N}a_{n}\left( \underset{v\in S_{n,V}}{\sum }\frac{\deg (v)}{%
\left| V(G)\right| }+\underset{w\in S_{n,FP}}{\sum }W(\delta
_{w}^{r})f(\delta _{w}^{r})\right) .
\end{array}
$
\end{center}

\strut

On the set of all simple functions, we can determine the ``almost
everywhere'' relation (in short, ``a.e'' relation). \strut

\strut

\begin{proposition}
Let $g_{1}$ $=\sum_{j=1}^{n}1_{S_{j}}$ and $g_{2}=\sum_{i=1}^{m}1_{T_{i}}$
be simple functions, where $S_{j}$'s and $T_{i}$'s are subsets of $\Bbb{F}%
^{+}(G^{\symbol{94}}).$ Suppose that $S_{j}$'s are mutually disjoint and
also $T_{i}$'s are mutually disjoint. If $\cup _{j=1}^{n}$ $S_{j}$ $=$ $\cup
_{i=1}^{m}$ $T_{i}$ in $\Bbb{F}^{+}(G^{\symbol{94}}),$ then $I(g_{1})$ $=$ $%
I(g_{2}).$
\end{proposition}

$\strut $

\begin{proof}
Observe that

\strut

$\ \ \ \ \ \ \ \ \ \ I(g_{1})=\sum_{j=1}^{n}\left( \underset{v\in S_{j,V}}{%
\sum }\frac{\deg (v)}{\left| V(G)\right| }+\underset{w\in S_{j,FP}}{\sum }%
W(\delta _{w}^{r})f(\delta _{w}^{r})\right) $

and

$\ \ \ \ \ \ \ \ \ \ I(g_{2})=\sum_{i=1}^{m}\left( \underset{v\in T_{i,V}}{%
\sum }\frac{\deg (v)}{\left| V(G)\right| }+\underset{w\in T_{i,FP}}{\sum }%
W(\delta _{w}^{r})f(\delta _{w}^{r})\right) .$

\strut \strut \strut

Since $\{S_{j}\}_{j=1}^{n}$ and $\{T_{i}\}_{i=1}^{m}$ are disjoint families
of $P\left( \Bbb{F}^{+}(G^{\symbol{94}})\right) ,$

\strut \strut

$\ \ \ \ \ \ \ \ \ \ \ \ \ \ \ \ \ \ \ \sum_{j=1}^{n}\underset{v\in S_{j,V}}{%
\sum }\frac{\deg (v)}{\left| V(G)\right| }=\underset{v\in \cup
_{j=1}^{n}S_{j,V}}{\sum }\frac{\deg (v)}{\left| V(G)\right| }$

and

$\ \ \ \ \ \ \ \ \ \ \ \ \ \ \ \ \ \ \ \sum_{i=1}^{m}\underset{v\in T_{i,V}}{%
\sum }\frac{\deg (v)}{\left| V(G)\right| }=\underset{v\in \cup
_{i=1}^{m}T_{i,V}}{\sum }\frac{\deg (v)}{\left| V(G)\right| },$

\strut

by the mutually disjointness on $\{S_{j,V}\}_{j=1}^{n}$ and $%
\{T_{i,V}\}_{i=1}^{m}.$ Similarly, by the mutually disjointness of $%
\{S_{j,FP}\}_{j=1}^{n}$ and $\{T_{i,FP}\}_{i=1}^{m}$ and by (3.7),

\strut

$\ \ \ \ \ \ \ \ \ \ \ \ \ \ \ \sum_{j=1}^{n}\underset{w\in S_{j,FP}}{\sum }%
\omega (\delta _{w}^{r})=\underset{w\in \cup _{j=1}^{n}S_{j,FP}}{\sum }%
\omega (\delta _{w}^{r})$

and

$\ \ \ \ \ \ \ \ \ \ \ \ \ \ \ \sum_{i=1}^{m}\underset{w\in T_{i,FP}}{\sum }%
\omega (\delta _{w}^{r})=\underset{w\in \cup _{i=1}^{m}T_{i,FP}}{\sum }%
\omega (\delta _{w}^{r}).$

\strut \strut

Therefore, we can conclude that $I(g_{1})=I(g_{2}),$ by (3.5).
\end{proof}

\strut

It is easy to check that, if $g_{j}$ $=$ $\sum_{k=1}^{n}a_{j,k}1_{S_{j,k}},$
for $j$ $=$ $1,$ $2,$ then we can have

\strut

(5.3) $\ \ \ \ \ I_{G}\left( g_{1}+g_{2}\right) =I_{G}(g_{1})+I_{G}(g_{2})$

(5.4) \ \ \ $\ \ I_{G}(c\cdot g_{j})=c\cdot I_{G}(g_{j}),$ for all $j=1,2.$

\strut

Now, let $1_{S}$ and $1_{T}$ be the characteristic functions, where $S$ $%
\neq $ $T$ in $P\left( \Bbb{F}^{+}(G^{\symbol{94}})\right) .$ Then we have
that $1_{S}$ $\cdot $ $1_{T}$ $=$ $1_{S\cap T}.$

\strut \strut

\begin{proposition}
Let $g_{j}$ be given as above, for $j$ $=$ $1,$ $2.$ Then

\strut 

$\ I_{G}(g_{1}g_{2})=\sum_{k,i=1}^{n}\left( a_{1,k}a_{2,i}\right) $

\strut 

$\ \ \ \ \ \ \ \ \ \ \ \ \ \ \ \ \ \ \cdot \left( \underset{v\in
S_{(1,k),V}\cap S_{(2,i),V}}{\sum }\frac{\deg (v)}{\left| V(G)\right| }+%
\underset{w\in S_{(1,k),FP}\cap S_{(2,i),FP}}{\sum }\omega (\delta
_{w}^{r})\right) .$
\end{proposition}

\strut

\begin{proof}
Observe that

\strut

$\ \ \ g_{1}g_{2}=\left( \sum_{k=1}^{n}a_{1,k}1_{S_{1,k}}\right) \left(
\sum_{i=1}^{n}a_{2,i}1_{S_{2,i}}\right) $

\strut

$\ \ \ \ \ \ \ \ \ \ \ =\sum_{k,i=1}^{n}a_{1,k}a_{2,i}\left(
1_{S_{1,k}}\cdot 1_{S_{2,i}}\right)
=\sum_{k,i=1}^{n}a_{1,k}a_{2,i}1_{S_{1,k}\cap S_{2,i}}$

\strut

So, we have that

\strut

(5.5) $\ \ \ \ \ \ \ \ \ I(g_{1}g_{2})=\sum_{k,i=1}^{n}a_{1,k}a_{2,i}$ $\
\mu \left( S_{1,k}\cap S_{2,i}\right) .$

\strut

Consider

\strut

$\ \ \mu \left( S_{1,k}\cap S_{2,i}\right) =d\left( (S_{1,k}\cap
S_{2,i})_{V}\right) +\omega \left( (S_{1,k}\cap S_{2,i})_{FP}\right) $

\strut

$\ \ \ \ \ \ \ \ \ \ \ =d\left( S_{(1,k),V}\cap S_{(2,i),V}\right) +\omega
\left( S_{(1,k),FP}\cap S_{(2,k),FP}\right) $

\strut

$\ \ \ \ \ \ \ \ \ \ \ =\underset{v\in S_{(1,k),V}\cap S_{(2,i),V}}{\sum }%
\frac{\deg (v)}{\left| V(G)\right| }+\underset{w\in S_{(1,k),FP}\cap
S_{(2,i),FP}}{\sum }\omega (\delta _{w}^{r}).$

$\strut $
\end{proof}

\strut

Suppose that $g_{1}$ and $g_{2}$ are given as above and assume that the
families $\{S_{1,k}\}_{k=1}^{n}$ and $\{S_{2,i}\}_{i=1}^{n}$ are disjoint
family. i.e.,

$\strut $

\begin{center}
$\left( \cup _{k=1}^{n}S_{1,k}\right) \cap \left( \cup
_{i=1}^{n}S_{2,i}\right) =\varnothing .$
\end{center}

\strut

Then $I(g_{1}g_{2})=0.$

\strut

Similar to Chapter 3, we can define a $\mu $-measurable function $g_{w},$
for all $w$ $\in $ $\Bbb{F}^{+}(G^{\symbol{94}}),$ by

\strut

\begin{center}
$g_{w}=1_{S_{l}^{w}\cup S_{r}^{w}},$
\end{center}

where

\begin{center}
$S_{l}^{w}=\{w^{\prime }\in \Bbb{F}^{+}(G^{\symbol{94}}):ww^{\prime }\in 
\Bbb{F}^{+}(G^{\symbol{94}})\}$
\end{center}

and

\begin{center}
$S_{r}^{w}=\{w^{\prime \prime }\in \Bbb{F}^{+}(G^{\symbol{94}}):w^{\prime
\prime }w\in \Bbb{F}^{+}(G^{\symbol{94}})\}.$
\end{center}

\strut \strut

\begin{proposition}
Let $w\in \Bbb{F}^{+}(G^{\symbol{94}})$ and assume that $g_{w}$ is the above 
$\mu $-measurable map. Then

\strut 

(5.6)$\ \ \ \ \ \ \ \ \ \ \ \ \ \ \ \ \ \ \ I(g_{w})$ $=$ $\mu $ $\left(
S_{l}^{w}\cup S_{r}^{w}\right) ,$

\strut 

where $S_{l}^{w}$ and $S_{r}^{w}$ are given in the previous paragraph. $%
\square $
\end{proposition}

\strut \strut 

So, different from Chapter 3, if $w$ is a loop finite path in $\Bbb{F}%
^{+}(G^{\symbol{94}}),$ then the integral $I(g_{w})$ of $w$ is $\infty ,$
whenever the weight $W(\delta _{w}^{r})$ is large enough in the interval $(0,
$ $1),$ because both $S_{l}^{w}$ and $S_{r}^{w}$ contains $w^{k},$ for all $k
$ $\in $ $\Bbb{N}$. 

\strut 

Let $g$ be a $\mu $-measurable function. Then the support $\Bbb{F}^{+}(G^{%
\symbol{94}}$ $:$ $g)$ of $g$ is well-defined. Define the monomial $g_{k}(w)$
$=$ $g_{w^{k}}$ on $\Bbb{F}^{+}(G^{\symbol{94}}),$ for $k$ $\in $ $\Bbb{N}.$
Then 

\strut 

\begin{center}
$\Bbb{F}^{+}(G^{\symbol{94}}:g_{1})=\Bbb{F}^{+}(G^{\symbol{94}})$
\end{center}

and

\begin{center}
$\Bbb{F}^{+}(G^{\symbol{94}}:g_{n})=V(G^{\symbol{94}})\cup loop(G^{\symbol{94%
}}),$
\end{center}

\strut \strut 

for all $n$ $\in $ $\Bbb{N}$ $\setminus $ $\{1\},$ where 

\begin{center}
$loop(G^{\symbol{94}})\overset{def}{=}\{l\in \Bbb{F}^{+}(G^{\symbol{94}}):$ $%
l$ is a loop finite path$\}.$
\end{center}

\strut 

So, by (5.6) and the supports, $I(g_{k})$ $=$ $\infty $ $=$ $I(g),$ for $k$ $%
\in $ $\Bbb{N},$ in general, where $g$ $=$ $\sum_{k=0}^{N}$ $a_{k}$ $g_{k}$
is a polynomial, with $g_{0}$ $\equiv $ $1.$  

\strut 

Define the $\mu $-measurable map $g_{-n}(w)$ $=$ $g_{w^{-n}},$ for all $w$ $%
\in $ $\Bbb{F}^{+}(G^{\symbol{94}})$ and $n$ $\in $ $\Bbb{N}.$ Then similar
to Chapter 3, we have that

\strut 

(5.7) $\ \ \ \ \ \ \ \ \ \ \ \Bbb{F}^{+}(G^{\symbol{94}}:g_{-n})=\Bbb{F}%
^{+}(G^{\symbol{94}}:g_{n}),$ for $n$ $\in $ $\Bbb{N}.$

\strut 

(5.8)  \ \ \ \ \ $\Bbb{F}^{+}(G^{\symbol{94}}:g_{\pm n})=\Bbb{F}^{+}(G^{%
\symbol{94}}:g_{2}),$ for all $n$ $\in $ $\Bbb{N}$ $\setminus $ $\{1\}.$

\strut 

By (5.7) and (5.8), we have the following proposition.

\strut 

\begin{proposition}
$I\left( g_{-n}\right) =I\left( g_{n}\right) ,$ for all $n$ $\in $ $\Bbb{N}.$
In particular, $I\left( g_{\pm n}\right) $ $=$ $I\left( g_{2}\right) ,$ for
all $n$ $\in $ $\Bbb{N}$ $\setminus $ $\{1\}.$ $\square $
\end{proposition}

\strut \strut \strut 

\strut \strut 

\strut \strut \strut \strut \strut 

\strut \textbf{References}

\strut

\strut

\begin{quote}
{\small [1] \ \ A. G. Myasnikov and V. Shapilrain (editors), Group Theory,
Statistics and Cryptography, Contemporary Math, 360, (2003) AMS.}

{\small [2] \ \ D. G. Radcliffe, Rigidity of Graph Products of Groups, Alg
\& Geom. Topology, Vol 3, (2003) 1079 - 1088.}

{\small [3] \ \ D. Voiculescu, Entropy of Random Walks on Groups and the
Macaev Norm, Proc. AMS, Vol 119, 3, (1993), 971 - 977.}

{\small [4] \ \ E. Breuillard and T. Gelander, Cheeger Constant and
Algebraic Entropy of Linear Groups, (2005) Preprint.}

{\small [5] \ \ F. Balacheff, Volum Entropy, Systole and Stable Norm on
Graphs, (2004) Preprint.}

{\small [6] \ \ G. C. Bell, Growth of the Asymptotic Dimension Function for
Groups, (2005) Preprint.}

{\small [7] \ \ I. Cho, Random Variables in a Graph }$W^{*}${\small %
-Probability Space, Ph. D. Thesis, (2005) Univ. of Iowa.}

{\small [8] \ \ I. Cho, Noncommutative Probability on Graph von Neumann
Algebras, (2006) Preprint.}

{\small [9] \ \ I. Cho, Compressed Random Variables in a Graph }$W^{*}$%
{\small -Probability Space, (2005), Preprint.}

{\small [10] J. Friedman and J-P. Tillich, Calculus on Graphs, (2005)
Preprint.}

{\small [11] J. Stallings, Centerless Groups-An Algebraic Formulation of
Gottlieb's Theorem, Topology, Vol 4, (1965) 129 - 134.}

{\small [12] R. G. Bartle, The Elements of Integration, (1966) 1-st Edition,
John Wiley \& Sons.}

{\small [13] R. Gliman, V. Shpilrain and A. G. Myasnikov (editors),
Computational and Statistical Group Theory, Contemporary Math, 298, (2001)
AMS.}
\end{quote}

\end{document}